\documentclass{article}

\usepackage{arxiv}

\usepackage[utf8]{inputenc} 
\usepackage[T1]{fontenc}    
\usepackage{hyperref}       
\usepackage{url}            
\usepackage{booktabs}       
\usepackage{amsfonts}       
\usepackage{nicefrac}       
\usepackage{microtype}      
\usepackage{graphicx}

\title{Analytical Approximations in Short Times of Exact Operational Solutions to Reaction–Diffusion Problems on Bounded Intervals}

\author{ \href{https://orcid.org/0000-0001-6679-4569}{\includegraphics[scale=0.06]{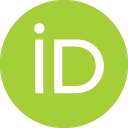}\hspace{1mm}Kwassi Anani}\thanks{ananikwassi@yahoo.fr} \\
	Department of Mathematics\\
	University of Lom\'e\\
	02 BP 1515, Maritime District, Lom\'e, Togo\\
	\texttt{kanani@univ-lome.tg} \\

}

\date{}

\renewcommand{\headeright}{}
\renewcommand{\undertitle}{}

\begin{document}
\maketitle

\begin{abstract}
This paper aims at obtaining, by means of integral transforms, analytical approximations in short times of solutions to boundary value problems for the one-dimensional reaction-diffusion equation with constant coefficients. The general form of the equation is considered on a bounded generic interval and the three classical types of boundary conditions, i.e., Dirichlet as well as Neumann and mixed boundary conditions are considered in a unified way. The Fourier and Laplace integral transforms are successively applied and an exact solution is obtained in the Laplace domain. This operational solution is proven to be the accurate Laplace transform of the infinite series obtained by the Fourier decomposition method and presented in the literature as solutions to this type of problem. On the basis of this unified operational solution, four cases are distinguished where innovative formulas expressing consistent analytical approximations in short time limits are derived with respect to the behavior of the solution at the boundaries. Compared to the infinite series solutions, the analytical approximations may open new perspectives and applications, among which can be noted the improvement of numerical efficiency in simulations of one-dimensional moving boundary problems, such as in Stefan models.
\end{abstract}


\section{Introduction }\label{intro}
Initial boundary value problems for linear parabolic equations are still much used nowadays as first model approximations of nonlinear and time-dependent problems in bounded domains (see for examples \cite{Br-analytical} and  \cite{Ca-Ma-semi}). Until now, exact analytical solutions to such linear boundary value problems are mainly obtained as infinite series of functions as in \cite{Lu-analytical} or in \cite{Hen-partial}. But, accurate analytical approximations in short times limits are hardly derivable from such infinite series solutions (see \cite{Ur-Pa-improved} and \cite{Ch-Ch-Ts-hermite}). The series solutions are generally obtained via the Fourier decomposition method, by means of separation of variables and using the regular Sturm-Liouville theory, as presented in many textbooks of which \cite{Ha-analytical} and \cite{Do-applied} can be cited among others. As shown in \cite{Hen-partial}, any one-dimensional linear parabolic equation with constant coefficients and a convection term can be reduced to an equation without convection term. Similarly, the linear reaction term with a fixed coefficient can also be canceled by specific transformations of the equation. Still, this latter reduction is not made here for reasons of later generalization. 

Compared to the classical numerical schemes as reported in \cite{Mi-handbook}, recent numerical approaches to boundary value problems for heat transfer equations have much  gained in sophistication and accuracy as it can be seen in \cite{Li-Xu-Zh-simulation}, \cite{Ta-Ya-Yo-Sr-computational} and \cite{Li-Wa-Wu-stable} for examples. However, apart of the Fourier decomposition method mentioned above, there exist no other established procedure for finding exact analytical solutions to boundary value problems involving linear reaction-diffusion equations even with constant coefficients. According to \cite{Lu-analytical}, the Laplace integral transform and its inversion formula are found to be not appropriate for solving boundary value problems with a non-uniform function as initial condition. Moreover, the classical Fourier integral applied to space coordinates is validated only for infinite and semi-infinite solids. Concerning other analytical and semi-approximate methods, some powerful techniques have been used to handle heat transfer problems in finite or infinite domains but without exhibiting new exact solutions. One can cite the Adomian’s decomposition method as in \cite{Mo-algorithm} and \cite{Bo-Ra-Wa-Du-adomian}, the variational iteration method as in \cite{Li-Zh-he}, the homotopy perturbation method as in \cite{Gh-Ka-applications} and Bessel collocation method as in \cite{Yu-Sa-numerical}. More details on analytical methods applicable to all types of equations are available in the book by \cite{Zh-modeling}.

In the present paper, an exact solution in the Laplace domain with computationally efficient analytical approximations in short time limits are obtained to the boundary value problem for the linear reaction-diffusion equation with constant coefficients. In section 2, the general form of the problem is introduced and the assumptions are briefly commented. Section 3 presents the application of the Fourier integral transform and the resulting solution is written under the form of integrals. In section 4, an explicit  solution in the Laplace domain is deduced from the previous integral form of the solution. Then, analytical solutions in short time limits are obtained in section 5, especially with respect to the behavior of the operational solution at the ends points of the domain. The results are briefly applied to a specific example in section 6 and the unified solution in the Laplace domain has been extended to unbounded space domains. Finally, section 7 outlines the conclusion.

\section{The reaction–diffusion equation}\label{sec:2}
We consider the one-dimensional linear parabolic equation with constant coefficients, expressed in the  nonhomogeneous reaction-diffusion form on a bounded generic interval. For space variable $ x \in [l_1, l_2]$, ($l_1, l_2 $ reals, $l_1< l_2$) , and for time $ t \in [0, T] $, ($T>0 $ or $t>0 $ in the case when $ T= +\infty$), the equation is stated as:
\begin{equation}
	\frac{\partial u}{\partial t}-a^2\frac{\partial^{2}u}{\partial x^{2}} + bu =f(x,t), \quad l_1 < x < l_2,\quad 0 < t < T.
	\label{eq01}
\end{equation}
The equation (\ref{eq01}) is subject to the initial condition
\begin{equation}
	u(x,0)=\varphi(x), \quad    l_1 \le x \le l_2,
	\label{eq02}
\end{equation}
and to the boundary conditions:
\begin{equation}
	\alpha_1 u(l_1,t)+\beta_1 \frac{\partial u}{\partial x}(l_1,t)=g_1(t), \quad \alpha_1, \beta_1 \in \mathbb{R}, \quad \alpha_1^2+\beta_1^2 \ne 0,
	\label{eq03}
\end{equation}
\begin{equation}
	\alpha_2 u(l_2,t)+\beta_2 \frac{\partial u}{\partial x}(l_2,t)=g_2(t),  \quad \alpha_2, \beta_2 \in \mathbb{R}, \quad \alpha_2^2+\beta_2^2 \ne 0,
	\label{eq04}
\end{equation}
for $ t \in [0, T] $. The four terms in equation (\ref{eq01}) represent respectively transient, diffusion, reaction and source terms. The function $ u(x,t) $ is to be determined, and may represent species concentration for mass transfer or temperature for heat transfer, while the functions $ f(x,t), g_1 (t)$ and $g_2(t) $ are given. The coefficient $ a>0 $ is related to the constant diffusivity of the mass or heat transfer. The reaction term (linear term in $u$) indicates the possibility of mass or heat exchange with the environment through the lateral surface of the body, at the rates proportional to the concentration or to the temperature \cite{He-mathematical}. In a process of mass diffusion for example, $ b $ is the coefficient of disintegration ($ b < 0 $) or multiplication ($ b > 0 $). Without lost of generality, we will assume from now on that $ b \geq 0 $. The source term expression $ f(x,t) $ may permit to consider homogeneous as well as nonhomogeneous equations, while the function $ \varphi(x) $ may indicate zero or non-zero initial condition. 

The Dirichlet, Neumann and mixed boundary conditions are expressed in a unified way by equations (\ref{eq03}) and (\ref{eq04}), that is, both homogeneous and nonhomogeneous forms of those three types of boundary conditions can be taken into account. It is sufficient to give some acceptable values to the real parameters $ \alpha_1, \beta_1, \alpha_2, \beta_2 $, and some convenient expressions to the time-dependent functions $ g_1 (t) $ and $g_2(t)$. Thus,  Dirichlet  conditions are satisfied when $\alpha_1=\alpha_2=1$ and  $\beta_1 =\beta_2=0 $, while Neumann conditions can be obtained on the boundaries if $\alpha_1=\alpha_2=0$ and  $\beta_1 =\beta_2=1 $. Likewise mixed boundary conditions are obtained when $\alpha_1=\alpha_2=1$ and  $\beta_1 \ne 0, \ \beta_2\ne 0 $. As it can be noted, a given combination of these three classical types of boundary conditions can also be obtained. Classically, the problem (\ref{eq01})-(\ref{eq04}) is first solved for homogeneous equation and boundary conditions ($ f(x,t)=0 $ and $ g_1 (t) =g_2(t)=0$), by using the Fourier decomposition method through the separation of variables, together with the Sturm-Liouville theory of eigenvalues and eigenfunctions. Then, the principle of Duhamel intervenes in addition, when accounting for a non-zero source term $ f(x,t) $. Now, in the case of nonhomogeneous boundary conditions, the problem (\ref{eq01})-(\ref{eq04}) is first reduced to a problem with boundary conditions equal to zero by the means of the so-called auxiliary functions (see \cite{Br-analytical} for example). 

The application of the above-mentioned techniques to the problem (\ref{eq01})-(\ref{eq04}) will lead to an exact series solution, which converges uniformly as well as the series obtained by differentiating twice by $ x $ and once by $ t  $, provided that, $ \varphi(x) $ and $ f(x,t) $  are assumed  continuous on $ [l_1, l_2]$ and on $[l_1, l_2]\times [0, T]$, and $ g_1 (t) $ and $g_2(t)$ are continuously differentiable on $ [0, T] $. The uniqueness of the series solution can be proven by using the maximum principle. In general, existence and uniqueness theorems of classical as well as weak or generalized solutions to initial-boundary-value problems for the linear parabolic equation in one or higher dimensional spaces, have been proven only under certain assumptions for specific classes of functions (see \cite{Ma-essential} and \cite{Ci-introduction} as  examples). Let us admit here in a first step that the source term $ f(x,t) $  is continuous on $[l_1, l_2]\times [0, T]$ and the initial condition function $ \varphi(x) $ satisfies at least the so-called Dirichlet conditions relatively to the space variable $ x $ on $ [l_1, l_2]$. That is, this function is piece-wise continuous or can be expressed in a unique way as a convergent series of eigenvalues and eigenfunctions of a Sturm-Liouville problem. Similarly, the time-dependent functions $ g_1 (t) $ and $g_2(t)$ are assumed here to be at least once piece-wise differentiable. Thus, they are continuous and their derivatives are piece-wise continuous on $ t\in [0, T] $. Standing on these hypotheses, a new method combining Fourier and Laplace integral transformations is proposed here, in order to obtain consistent analytical approximations in short time limits from the exact operational solution to the boundary value problem (\ref{eq01})-(\ref{eq04}). This implies that all the involved functions are assumed here to be absolutely integrable with respect to time and space variables, while the source term $ f(x,t) $ is assumed to be of exponential order relatively to the time variable $ t $, as well as the functions $ g_1 (t) $ and $g_2(t)$ and their respective derivatives.

\section{Method of Fourier integral transform}\label{sec:3}
In this section, in order to obtain the solution to the problem (\ref{eq01})-(\ref{eq04}) under a form of an integral, the Fourier integral transform  (abbreviated FIT from now on) and its inverse are used in relation to the space variable $x$.
We first recall that for any absolutely integrable function $ \phi(x) $, i.e., $ \int^{+\infty}_{-\infty}\vert \phi(x)\vert dx< \infty $, the FIT $ \Phi $ is defined as:
\begin{equation} \Phi(\lambda)=\frac{1}{\sqrt{2\pi}}\int^{+\infty}_{-\infty}\phi(x)\exp(-i\lambda x)dx,
	\label{eq05}
\end{equation}
where $ i^2=-1 $, $ \lambda \in \mathbb{R}$ and $ x\in \mathbb{R} $. The inverse Fourier transform of $ \Phi $ is:
\begin{equation} \phi(x)=\frac{1}{\sqrt{2\pi}}\int^{+\infty}_{-\infty}\Phi(\lambda)\exp(ix\lambda) d\lambda.
	\label{eq06}
\end{equation}
As a basic property, the FIT  tends to 0 when $\vert\lambda\vert$ goes to $\infty$.

The homogeneous form of equation (\ref{eq01}) reads:
\begin{equation}
	\frac{\partial u}{\partial t}=a^2\frac{\partial^{2}u}{\partial x^{2}} - bu. 
	\label{eq07}
\end{equation}
Let $ u_1 (x,t) $ be a solution to the boundary value problem formed by the homogeneous equation (\ref{eq07}) together with the initial condition (\ref{eq02}), and the boundary conditions (\ref{eq03}) and (\ref{eq04}). Assuming that the functions  $ u_1 (x,t) $, $ 	\frac{\partial u_1}{\partial t}(x,t) $ and $ \frac{\partial^{2}u_1}{\partial x^{2}}(x,t) $ are absolutely integrable with respect to the variables $x$ and $ t $, they  can be identified to their extension by $ 0 $ outside the rectangle $[l_1, l_2]\times [0, T] $, without lost of generality. The Fourier integral transform (FIT) with respect to the space variable $ x $ can be applied to $ u_1 (x,t) $ and will give:
$$F(\lambda,t)=\frac{1}{\sqrt{2\pi}}\int^{+\infty}_{-\infty}u_1 (x,t) \exp(-i\lambda x)dx =\frac{1}{\sqrt{2\pi}}\int^{l_2}_{l_1}u_1 (x,t)\exp(-i\lambda x)dx.$$
The FIT applied to the transient term $ 	\frac{\partial u_1}{\partial t}(x,t) $ leads to:
$$A(\lambda,t)=\frac{1}{\sqrt{2\pi}}\int^{l_2}_{l_1}\frac{\partial u_1}{\partial t}(x,t)\exp(-i\lambda x)dx=\frac{1}{\sqrt{2\pi}}\frac{\partial}{\partial t}\int^{l_2}_{l_1}u_1(x,t)\exp(-i\lambda x)dx=\frac{\partial F}{\partial t}(\lambda,t)$$
Similarly, the FIT of the diffusion term $ a^2\frac{\partial^{2}u_1}{\partial x^{2}}(x,t) $ can be written as:
$$B(\lambda,t)=\frac{a^2}{\sqrt{2\pi}}\int^{l_2}_{l_1}\frac{\partial^{2}u_1}{\partial x^{2}}(x,t)\exp(-i\lambda x)dx.$$
By using two successive integration by parts over the finite interval $[l_1, l_2]$, $ B(\lambda,t) $ can be expressed as:
\begin{equation}
	\begin{array}{ll}
		\displaystyle B(\lambda,t)=-a^2 \lambda^2 F(\lambda,t)+\frac{a^2}{\sqrt{2\pi}}
		\left[\frac{\partial u_1}{\partial x}(l_2,t)\exp(-i\lambda l_2)
		-\frac{\partial u_1}{\partial x}(l_1,t)\exp(-i\lambda l_1)\right]\\
		\displaystyle+\frac{a^2}{\sqrt{2\pi}}
		\left[i\lambda u_1(l_2,t)\exp(-i\lambda l_2)-i\lambda u_1(l_1,t)\exp(-i\lambda l_1)\right].
	\end{array}
	\label{eq08}
\end{equation}
The FIT of the linear term $ -bu_1(x,t)  $ is simply
$$C(\lambda,t)=-b F(\lambda,t).$$
Now, according to the application of the FIT to the terms of the homogeneous equation, if $ u_1 (x,t)$, considered to be null outside $[l_1, l_2]\times [0, T]  $ is solution of equation (\ref{eq07}), then $ F(\lambda,t) $ is solution of the equation: $ A(\lambda,t)- B(\lambda,t)-C(\lambda,t)=0$,  i.e., 
\begin{equation}
	\begin{array}{ll}
		\displaystyle \frac{\partial F}{\partial t}(\lambda,t)+(b+a^2 \lambda^2) F(\lambda,t)=\frac{a^2}{\sqrt{2\pi}}
		\left[\frac{\partial u_1}{\partial x}(l_2,t)\exp(-i\lambda l_2)
		-\frac{\partial u_1}{\partial x}(l_1,t)\exp(-i\lambda l_1)\right]\\
		\displaystyle+\frac{a^2}{\sqrt{2\pi}}
		\left[i\lambda u_1(l_2,t)\exp(-i\lambda l_2)-i\lambda u_1(l_1,t)\exp(-i\lambda l_1)\right];
	\end{array}
	\label{eq09}
\end{equation}
for $ \lambda \in \mathbb{R} $ and $ 0 \le t \le T $.

In order to integrate equation (\ref{eq09}), we multiply each member by $ \exp [(b+a^2 \lambda^2)t] $ and obtain:
\begin{equation}
	\begin{array}{ll}
		\displaystyle \frac{\partial }{\partial t}(F(\lambda,t)\exp [(b+a^2 \lambda^2)t])=\\
		\displaystyle\frac{a^2}{\sqrt{2\pi}}
		\exp [(b+a^2 \lambda^2)t]\left[\frac{\partial u_1}{\partial x}(l_2,t)\exp(-i\lambda l_2)
		-\frac{\partial u_1}{\partial x}(l_1,t)\exp(-i\lambda l_1)\right]\\
		\displaystyle+\frac{a^2}{\sqrt{2\pi}}
		\exp [(b+a^2 \lambda^2)t]\left[i\lambda u_1(l_2,t)\exp(-i\lambda l_2)-i\lambda u_1(l_1,t)\exp(-i\lambda l_1)\right].
	\end{array}
	\label{eq10}
\end{equation} 
Proceeding now by integration of equation (\ref{eq10}) relatively to the time variable from $ \eta=0 $ to $ \eta=t \le T $, we obtain:
\begin{equation}
	\begin{array}{ll}
		\displaystyle F(\lambda,t)\exp [(b+a^2 \lambda^2)t]-F(\lambda,0)=\\
		\displaystyle \frac{a^2}{\sqrt{2\pi}}
		\int_{0}^{t}\left[\frac{\partial u_1}{\partial x}(l_2,\eta)\exp(-i\lambda l_2)
		-\frac{\partial u_1}{\partial x}(l_1,\eta)\exp(-i\lambda l_1)\right]\exp [(b+a^2 \lambda^2)\eta]d\eta\\
		\displaystyle+\frac{a^2}{\sqrt{2\pi}}
		\int_{0}^{t}\left[i\lambda u_1(l_2,\eta)\exp(-i\lambda l_2)-i\lambda u_1(l_1,\eta)\exp(-i\lambda l_1)\right]\exp [(b+a^2 \lambda^2)\eta]d\eta.
	\end{array}
	\label{eq11}
\end{equation} 
Due to the initial condition (\ref{eq02}), $ F(\lambda,0)= F(\lambda,t=0)$ can be calculated as 
\begin{equation} F(\lambda,0)=\frac{1}{\sqrt{2\pi}}\int^{l_2}_{l_1}\varphi(\xi)\exp(-i\lambda \xi)d\xi,
	\label{eq12}
\end{equation}
where the dummy variable is replaced by the  integrating variable $\xi$ in order to avoid confusion.
Equation (\ref{eq11}) can then be rewritten as:
\begin{equation}
	\begin{array}{ll}
		\displaystyle F(\lambda,t)=\\
		\displaystyle \frac{a^2}{\sqrt{2\pi}}
		\int_{0}^{t}\left[\frac{\partial u_1}{\partial x}(l_2,\eta)\exp(-i\lambda l_2)
		-\frac{\partial u_1}{\partial x}(l_1,\eta)\exp(-i\lambda l_1)\right]\exp [-(b+a^2 \lambda^2)(t-\eta)]d\eta\\
		\displaystyle+\frac{a^2}{\sqrt{2\pi}}
		\int_{0}^{t}\left[i\lambda u_1(l_2,\eta)\exp(-i\lambda l_2)-i\lambda u_1(l_1,\eta)\exp(-i\lambda l_1)\right]\exp [-(b+a^2 \lambda^2)(t-\eta)]d\eta\\
		\displaystyle+\frac{1}{\sqrt{2\pi}}\exp [-(b+a^2 \lambda^2)t]\int^{l_2}_{l_1}\varphi(\xi)\exp(-i\lambda \xi)d\xi.
	\end{array}
	\label{eq13}
\end{equation} 

In order to obtain $ u_1(x,t) $, the inversion formula (\ref{eq06}) needs to be applied to the function $ F(\lambda,t) $ expressed by equation (\ref{eq13}), the Fourier variable $ \lambda $ running from $ -\infty $ to $ +\infty $. First, due to properties of the convolution of two functions, the three terms of the second member can be respectively denoted as:
$$ F_1(\lambda,t)=\frac{a^2}{\sqrt{2\pi}}
\int_{0}^{t}\left[\frac{\partial u_1}{\partial x}(l_2,t-\eta)\exp(-i\lambda l_2)
-\frac{\partial u_1}{\partial x}(l_1,t-\eta)\exp(-i\lambda l_1)\right]\exp [-(b+a^2 \lambda^2)\eta]d\eta, $$
$$ F_2(\lambda,t)=\frac{a^2}{\sqrt{2\pi}}
\int_{0}^{t}\left[i\lambda u_1(l_2,t-\eta)\exp(-i\lambda l_2)-i\lambda u_1(l_1,t-\eta)\exp(-i\lambda l_1)\right]\exp [-(b+a^2 \lambda^2)\eta]d\eta$$
and
$$ F_3(\lambda,t)=\frac{1}{\sqrt{2\pi}}\exp [-(b+a^2 \lambda^2)t]\int^{l_2}_{l_1}\varphi(\xi)\exp(-i\lambda \xi)d\xi. $$
Changing the order of integration due to the  convergence of the integrals involved, the inverse $ I_1 $ of the first term $ F_1 $ is calculated as:
$$ 
\begin{array}{ll}
	I_1(x,t)=\displaystyle\frac{a^2}{{2\pi}}
	\int_{0}^{t}\frac{\partial u_1}{\partial x}(l_2,t-\eta)\int_{ -\infty }^{\infty }\exp(-i\lambda l_2)\exp(i\lambda x)\exp [-(b+a^2 \lambda^2)\eta]d\lambda d\eta\\
	-\displaystyle\frac{a^2}{{2\pi}}\int_{0}^{t}\frac{\partial u_1}{\partial x}(l_1,t-\eta)\int_{ -\infty }^{\infty }\exp(-i\lambda l_1)\exp(i\lambda x)\exp [-(b+a^2 \lambda^2)\eta]d\lambda d\eta.  
\end{array}
$$
By means of computations $ I_1 $ is reduced to:
\begin{equation}
	\begin{array}{ll}
		I_1(x,t)=\displaystyle\frac{a}{2\sqrt{\pi}}
		\int_{0}^{t}\frac{\partial u_1}{\partial x}(l_2,t-\eta)\left[ \exp\left( -\frac{(l_2-x)^2}{4 a^2 \eta}\right) \right]\frac{\exp(-b\eta)}{\sqrt{\eta}} d\eta\\
		-\displaystyle\frac{a}{2\sqrt{\pi}}
		\int_{0}^{t}\frac{\partial u_1}{\partial x}(l_1,t-\eta)\left[ \exp\left(-\frac{(l_1-x)^2}{4 a^2 \eta}\right) \right]\frac{\exp(-b\eta)}{\sqrt{\eta}} d\eta.  
	\end{array}
	\label{eq14} 
\end{equation}
Similarly, $ I_2 $ is obtained as:
$$ \begin{array}{ll}
	I_2(x,t)=\displaystyle\frac{a^2}{{2\pi}}
	\int_{0}^{t}u_1(l_2,t-\eta)\int_{ -\infty }^{\infty }i\lambda\exp(-i\lambda l_2)\exp(i\lambda x)\exp [-(b+a^2 \lambda^2)\eta]d\lambda d\eta\\
	-\displaystyle\frac{a^2}{{2\pi}}\int_{0}^{t}u_1(l_1,t-\eta)\int_{ -\infty }^{\infty }i\lambda\exp(-i\lambda l_1)\exp(i\lambda x)\exp [-(b+a^2 \lambda^2)\eta]d\lambda d\eta,  
\end{array}
$$
i.e.,
\begin{equation}
	\begin{array}{ll} 
		I_2(x,t)=\displaystyle\frac{1}{4 a\sqrt{\pi}}
		\int_{0}^{t}u_1(l_2,t-\eta)\left[(l_2-x) \exp\left( -\frac{(l_2-x)^2}{4 a^2 \eta}\right) \right]\frac{\exp(-b\eta)}{{\eta}^{3/2}}  d\eta\\
		-\displaystyle\frac{1}{4 a\sqrt{\pi}}
		\int_{0}^{t}u_1(l_1,t-\eta)\left[(l_1-x)\exp\left(-\frac{(l_1-x)^2}{4 a^2 \eta}\right) \right]\frac{\exp(-b\eta)}{{\eta}^{3/2}}  d\eta.
	\end{array} 
	\label{eq15} 
\end{equation}
And, the inverse $ I_3 $ of $ F_3 $ is calculated as:
\begin{equation}
	I_3(x,t)=\displaystyle\frac{\exp(-bt)}{2 a\sqrt{\pi t}}
	\int_{l_1}^{l_2}\varphi(\xi)\left[\exp\left( -\frac{(\xi-x)^2}{4 a^2 t}\right) \right] d\xi.
	\label{eq16} 
\end{equation}

From the calculations above, the expression of $ u_1 $ is deduced as $  u_1 (x,t)=I_1(x,t)+I_2(x,t)+I_3(x,t) $, that is:
\begin{equation} 
	\begin{array}{ll} 
		u_1 (x,t)=\displaystyle\frac{a}{2\sqrt{\pi}}
		\int_{0}^{t}\frac{\partial u_1}{\partial x}(l_2,t-\eta)\left[ \exp\left( -\frac{(l_2-x)^2}{4 a^2 \eta}\right) \right]\frac{\exp(-b\eta)}{\sqrt{\eta}} d\eta\\
		-\displaystyle\frac{a}{2\sqrt{\pi}}
		\int_{0}^{t}\frac{\partial u_1}{\partial x}(l_1,t-\eta)\left[ \exp\left(-\frac{(l_1-x)^2}{4 a^2 \eta}\right) \right]\frac{\exp(-b\eta)}{\sqrt{\eta}} d\eta\\
		+\displaystyle\frac{1}{4 a\sqrt{\pi}}
		\int_{0}^{t}u_1(l_2,t-\eta)\left[(l_2-x) \exp\left( -\frac{(l_2-x)^2}{4 a^2 \eta}\right) \right]\frac{\exp(-b\eta)}{{\eta}^{3/2}}  d\eta\\
		-\displaystyle\frac{1}{4 a\sqrt{\pi}}
		\int_{0}^{t}u_1(l_1,t-\eta)\left[(l_1-x)\exp\left(-\frac{(l_1-x)^2}{4 a^2 \eta}\right) \right]\frac{\exp(-b\eta)}{{\eta}^{3/2}}  d\eta\\
		+\displaystyle\frac{\exp(-bt)}{2 a\sqrt{\pi t}}
		\int_{l_1}^{l_2}\varphi(\xi)\left[\exp\left( -\frac{(\xi-x)^2}{4 a^2 t}\right) \right] d\xi.
	\end{array} 
	\label{eq17}	
\end{equation} 
It can be verified that $ u_1 $, as expressed by the integral form (\ref{eq17}), is a solution to the homogeneous equation (\ref{eq07}). Moreover, the initial condition (\ref{eq02}) is satisfied by the solution (\ref{eq17}) since, when $ t\rightarrow 0 $, $ I_1(x,t) $ and $ I_2(x,t)$ vanish and $ \displaystyle\lim_{t\rightarrow 0} u_1(x,t)$ reduces to $ \displaystyle\lim_{t\rightarrow 0} I_{3}(x,t)$. Now, if we put
$$ G(x,\xi,t)= \frac{\exp(-bt)}{2 a\sqrt{\pi t}}\exp\left( -\frac{(\xi-x)^2}{4 a^2 t}\right),$$ 
then, due to the property of the normal or Gauss probability density:
$$\displaystyle\lim_{t\rightarrow 0}\int_{-\infty}^{\infty} G(x,\xi,t) d\xi=\displaystyle\lim_{t\rightarrow 0}\exp(-bt)=1,$$
and $ G(x,\xi,t)\rightarrow 0 $ as $t\rightarrow 0 $ at all points $ (x,\xi) \in \mathbb{R}^2 $, with the exception of the diagonal $ x=\xi $ where it becomes infinitely large. Therefore, if $ G(x,\xi,t) $ is an analogue of a Green's function and 
$$ \displaystyle \lim_{t\rightarrow 0} G(x,\xi,t)=\delta (x-\xi),$$ 
where $ \delta (x-\xi) $ is the Dirac delta function. If $ \mathbb{I}_{[l_1, l_2]} $ denotes the indicator function of interval $ [l_1, l_2] $, one has:
\begin{equation}
	\begin{array}{ll} 
		\displaystyle\lim_{t\rightarrow 0} u_1(x,t)= \displaystyle\lim_{t\rightarrow 0} I_{3}(x,t)=\displaystyle\lim_{t\rightarrow 0}\int_{l_1}^{l_2}\varphi(\xi)G(x,\xi,t)d\xi\\=\displaystyle\lim_{t\rightarrow 0}\int_{-\infty}^{\infty} \mathbb{I}_{[l_1, l_2]}(\xi)\varphi(\xi)G(x,\xi,t)d\xi
		=\int_{-\infty}^{\infty} \mathbb{I}_{[l_1, l_2]}(\xi)\varphi(\xi) \delta (x-\xi) d\xi=\varphi(x),
	\end{array} 
	\label{eq18}
\end{equation}
and the initial condition (\ref{eq02}) is satisfied by $ u_1 $. In fact, it has just been proven that the function $ I_3 $ satisfied the non-zero initial condition (\ref{eq02}). Moreover, it can also be verified that $ I_3 $ is a solution to the homogeneous equation (\ref{eq07}). According to the Duhamel's principle, a solution of the nonhomogeneous equation (\ref{eq01}) with zero initial condition can be written as:
\[  u_2(x,t)=\int_{0}^{t}d\theta\int_{l_1}^{l_2}\,G(x,\xi,t-\theta)f(\xi,\theta)d\xi.
\]

Finally, by the superposition principle, a solution $ u $ of the nonhomogeneous equation (\ref{eq01}) together with the non-zero initial condition (\ref{eq02}) can be expressed as $ u=u_1+u_2 $, namely:
\begin{equation} 
	\begin{array}{ll} 
		u (x,t)=\displaystyle\frac{a}{2\sqrt{\pi}}
		\int_{0}^{t}\frac{\partial u}{\partial x}(l_2,t-\eta)\left[ \exp\left( -\frac{(l_2-x)^2}{4 a^2 \eta}\right) \right]\frac{\exp(-b\eta)}{\sqrt{\eta}} d\eta\\
		-\displaystyle\frac{a}{2\sqrt{\pi}}
		\int_{0}^{t}\frac{\partial u}{\partial x}(l_1,t-\eta)\left[ \exp\left(-\frac{(l_1-x)^2}{4 a^2 \eta}\right) \right]\frac{\exp(-b\eta)}{\sqrt{\eta}} d\eta\\
		+\displaystyle\frac{1}{4 a\sqrt{\pi}}
		\int_{0}^{t}u(l_2,t-\eta)\left[(l_2-x) \exp\left( -\frac{(l_2-x)^2}{4 a^2 \eta}\right) \right]\frac{\exp(-b\eta)}{{\eta}^{3/2}}  d\eta\\
		-\displaystyle\frac{1}{4 a\sqrt{\pi}}
		\int_{0}^{t}u(l_1,t-\eta)\left[(l_1-x)\exp\left(-\frac{(l_1-x)^2}{4 a^2 \eta}\right) \right]\frac{\exp(-b\eta)}{{\eta}^{3/2}}  d\eta\\
		+r (x,t)
	\end{array} 
	\label{eq19}	
\end{equation} 
where:
\begin{equation} 
	\begin{array}{ll} 
		r (x,t)= \displaystyle\frac{\exp(-bt)}{2 a\sqrt{\pi t}}
		\int_{l_1}^{l_2}\varphi(\xi)\left[\exp\left( -\frac{(\xi-x)^2}{4 a^2 t}\right) \right] d\xi\\
		+\displaystyle\frac{1}{2 a\sqrt{\pi }}\int_{0}^{t}d\theta\int_{l_1}^{l_2}\frac{\exp(-b(t-\theta))}{\sqrt{(t-\theta)}}\exp\left( -\frac{(\xi-x)^2}{4 a^2 (t-\theta)}\right)f(\xi,\theta)d\xi
	\end{array} 
	\label{eq20}	
\end{equation} 
The FIT method has allowed to determine an expression in integral form of the solution of the equation (\ref{eq01}) under the initial condition (\ref{eq02}).  But, the solution (\ref{eq19}) depends of the values of the function $u$ and of its derivative at the boundaries. The boundary conditions (\ref{eq03}) and (\ref{eq04}) will be taken into account via the Laplace integral transform (abbreviated LIT).

\section{Exact solution in the Laplace domain}
\label{sec:4}
If $f (t)$ is a function defined in $t \geq 0$, then its unilateral Laplace integral transform (LIT) is given in the complex $ p $-plane by (see \cite{Her-partial}):
\begin{equation}
	F(p)={\mathcal L}\{f(t)\}=\int_{0}^\infty f(t)\mathrm{e}^{-pt} dt,
	\label{eq21}
\end{equation}
provided that $f (t)$ be of exponential order, that is, there are constants $C$ and $\sigma$ so that $|f(t)|<C \mathrm{e}^{\sigma t}$, when t is sufficiently large.
The inversion, from the Laplace domain $p$ to the time domain $t$ is given by the complex
integral,
\begin{equation}
	f(t)={\mathcal L}^{-1}\{F(p)\}=\frac{1}{2\pi i}\int_{\gamma-i\infty}^{\gamma+i\infty}F(p) \mathrm{e}^{pt} dp,
	\label{eq22}
\end{equation}
where $ \gamma > \sigma $ is chosen so that $ F(p) $ converges absolutely on the real part of $p$ line $ \Re(p) = \gamma $, and $F(p)$ is analytic at the right of this line. Also, the LIT and its inverse can be obtained for many usual functions by using tables of transforms as in \cite{Po-laplace}. In the case when $ p $ is real as considered in the sequel of this paper, the inequality  $ p\ge \sigma $ needs to be satisfied. Relatively to the LIT, an important property is the convolution theorem (see \cite{De-integral}):\\
Let $f (t)$ and $g (t)$ be functions defined in $t \geq 0$. If ${\mathcal L}\{f(t)\}= F(p)$ and ${\mathcal L}\{g(t)\}= G(p)$, then
$$
{\mathcal L}\{f(t)\ast g(t)\}= {\mathcal L}\{f(t)\}{\mathcal L}\{g(t)\}=F(p)G(p)
$$
where $f(t)\ast g(t)$ is called the convolution of $f(t)$ and $g(t)$ and is defined by
the integral
$$
f(t)\ast g(t)=\int^{t}_{0} f(t-\eta)g(\eta)d\eta.
$$

In order to apply the LIT, all the time dependent functions involved in the problem (\ref{eq01})-(\ref{eq04}) are assumed to be originals, that is, the transforms of these functions exist. Thus, the LIT of temperature distribution $ u(x,t) $,  source term $ f(x,t) $, boundary functions $ g_1(t) $ and $ g_2(t) $ are respectively denoted in the Laplace domain by $U(x,p)$, $F(x,p)$, $G_1(p)$ and $G_2(p)$. Under these hypotheses, the analogue of the initial boundary value problem (\ref{eq01})-(\ref{eq04}) in the Laplace domain can be written in the form of an ordinary differential equation as follows:
\begin{equation}
	-a^2\frac{d^2 U}{dx^{2}}(x,p)+(b+p)U(x,p)=F(x,p)+\varphi(x),
	\label{eq23}
\end{equation}
since the Laplace transform of the time derivative is equal to the product of the transform by the operator p minus the value of
the function at the initial time instant:
\begin{equation}
	{\mathcal L}\left\lbrace \frac{\partial u}{\partial t}(x,t)\right\rbrace= pU(x,p)-\varphi(x).
	\label{eq24}
\end{equation}
Equation (\ref{eq23}) is subjected to the boundary conditions:
\begin{equation}
	\alpha_1 U(l_1,p)+\beta_1 \frac{d U}{d x}(l_1,p)=G_1(p), \quad \alpha_1, \beta_1 \in \mathbb{R}, \quad \alpha_1^2+\beta_1^2 \ne 0,
	\label{eq25}
\end{equation}
\begin{equation}
	\alpha_2 U(l_2,p)+\beta_2 \frac{d U}{d x}(l_2,p)=G_2(p),  \quad \alpha_2, \beta_2 \in \mathbb{R}, \quad \alpha_2^2+\beta_2^2 \ne 0.
	\label{eq26}
\end{equation}

Due to the convolution theorem, the solution expressed by equation (\ref{eq19}) can also be transformed by the LIT to be written as:
\begin{equation}
	\begin{array}{ll} 
		U(x,p)=\\
		\displaystyle \frac{a}{2 \sqrt{b+p}}\left[U_x(l_2,p)\exp\left(\frac{-(l_2-x)\sqrt{b+p}}{a} \right)-U_x(l_1,p)\exp\left(\frac{-(x-l_1)\sqrt{b+p}}{a} \right)   \right]\\ 
		+\displaystyle \frac{1}{2}\left[U(l_2,p)\exp\left(\frac{-(l_2-x)\sqrt{b+p}}{a} \right)+U(l_1,p)\exp\left(\frac{-(x-l_1)\sqrt{b+p}}{a} \right)   \right]\\
		+\displaystyle R(x,p),
		\label{eq27}
	\end{array} 
\end{equation}
where $ U_x $ denotes the derivative $ \frac{dU}{dx} $ and $ R(x,p) $ stands for the Laplace transform of the remaining term $ r (x,t) $ expressed by the equation (\ref{eq20}), namely, $R(x,p)={\mathcal L}\{r(x,t)\}  $. Substituting respectively $ x=l_1 $ and then $ x=l_2 $ in the above expression of $ U(x,p) $, we have the following equations:
\begin{equation}
	\begin{array}{ll} 
		\displaystyle\frac{1}{2}U(l_1,p) + \frac{a}{2 \sqrt{b+p}}U_x(l_1,p) -\frac{1}{2}\exp\left(\frac{-(l_2-l_1)\sqrt{b+p}}{a} \right)U(l_2,p)\\
		\displaystyle - \frac{a}{2 \sqrt{b+p}}\exp\left(\frac{-(l_2-l_1)\sqrt{b+p}}{a} \right)U_x(l_2,p)= R(l_1,p),
		\label{eq28}
	\end{array} 
\end{equation}
and 
\begin{equation}
	\begin{array}{ll} 
		\displaystyle-\frac{1}{2}\exp\left(\frac{-(l_2-l_1)\sqrt{b+p}}{a} \right) U(l_1,p) + \frac{a}{2 \sqrt{b+p}}\exp\left(\frac{-(l_2-l_1)\sqrt{b+p}}{a} \right)U_x(l_1,p) \\
		\displaystyle+\frac{1}{2}U(l_2,p)  - \frac{a}{2 \sqrt{b+p}}  U_x(l_2,p)= R(l_2,p),
	\end{array} 
	\label{eq29}
\end{equation}
Now, equations (\ref{eq25}), (\ref{eq26}),  (\ref{eq28}) and (\ref{eq29}) form a system ($ \mathcal{S} $) of four linear equations with four unknown functions that are $U(l_1,p)$, $U_x(l_1,p)$, $U(l_2,p)  $  and $U(l_2,p)  $. Therefore, the function $ U(x,p) $ is solution of the analogue problem (\ref{eq23})-(\ref{eq26}) in the Laplace domain, or equivalently $ u(x,t) $ is solution of the problem (\ref{eq01})-(\ref{eq04}) in the time domain, if and only if, the system ($ \mathcal{S} $) admits a unique solution. The determinant of the system is calculated as:
\begin{equation}
	\mathrm{det}( \mathcal{S})=\left| 
	\begin{array}{cccc}
		\displaystyle \alpha_1	& \displaystyle\beta_1 & \displaystyle 0 &  \displaystyle 0\\
		\displaystyle 0	& 0 &\displaystyle \alpha_2 &  \displaystyle\beta_2\\
		\displaystyle\frac{1}{2}	& \displaystyle \frac{a}{2 \sqrt{b+p}} &  \displaystyle\frac{-\chi(l_2-l_1,p)}{2} & \displaystyle \frac{-a\chi(l_2-l_1,p)}{2 \sqrt{b+p}} \\
		\displaystyle\frac{-\chi(l_2-l_1,p)}{2}	& \displaystyle\frac{a\chi(l_2-l_1,p)}{2 \sqrt{b+p}} & \displaystyle\frac{1}{2} & \displaystyle \frac{-a}{2 \sqrt{b+p}}
	\end{array}
	\right|;  
	\label{eq30}
\end{equation}
where the function $\chi(x,p)  $ is defined as $ \chi(x,p)=\exp\left(\frac{-x\sqrt{b+p}}{a} \right) $. This determinant is reduced by computations into:
\begin{equation}
	\begin{array}{cccc}
		\displaystyle\det( \mathcal{S})= -\frac{1}{4}(a^2 \alpha_1 \alpha_2-a \alpha_1 \beta_2+a \alpha_2 \beta_1)\frac{\chi(2(l_2-l_1),p)}{b+p}+\frac{1}{4} \beta_1 \beta_2 \chi(2(l_2-l_1),p)\\
		\displaystyle+\frac{1}{4}\frac{a^2 \alpha_1 \alpha_2+a \alpha_1 \beta_2-a \alpha_2 \beta_1}{b+p}-\frac{1}{4} \beta_1 \beta_2.
	\end{array} 
	\label{eq31}
\end{equation}
Thus, the determinant of the system is null if and only if all the coefficients of the functions appearing in equation (\ref{eq31}) are null. This will imply:
\begin{equation}
	\left\lbrace 
	\begin{array}{cccc}
		a^2 \alpha_1 \alpha_2-a \alpha_1 \beta_2+a \alpha_2 \beta_1=0,\\
		a^2 \alpha_1 \alpha_2+a \alpha_1 \beta_2-a \alpha_2 \beta_1=0,\\
		\beta_1 \beta_2=0.
	\end{array}
	\right. 
	\label{eq32}
\end{equation}
The system (\ref{eq32}) can be split into two systems since it is equivalent to:
\begin{equation}
	\left\lbrace 
	\begin{array}{cccc}
		a^2 \alpha_1 \alpha_2-a \alpha_1 \beta_2=0,\\
		a^2 \alpha_1 \alpha_2+a \alpha_1 \beta_2=0,\\
		\beta_1=0,
	\end{array}
	\right. 
	\label{eq33}
\end{equation}
or 
\begin{equation}
	\left\lbrace 
	\begin{array}{cccc}
		a^2 \alpha_1 \alpha_2+a \alpha_2 \beta_1=0,\\
		a^2 \alpha_1 \alpha_2-a \alpha_2 \beta_1=0,\\
		\beta_2=0.
	\end{array}
	\right. 
	\label{eq34}
\end{equation}
Since $ a \ne 0 $, by adding and then subtracting the two first equations, the sub-system (\ref{eq33}) can be shown to be equivalent to:
\begin{equation}
	\left\lbrace 
	\begin{array}{cccc}
		\alpha_1 \alpha_2=0,\\
		\alpha_1 \beta_2=0,\\
		\beta_1=0.
	\end{array}
	\right. 
	\label{eq35}
\end{equation}
The system (\ref{eq35}) leads to a contradiction with the hypotheses on the coefficients $ \alpha_1, \:\beta_1,\:\alpha_2$ and $\beta_2$, since it implies that $ \alpha_1 =\beta_1=0$ or $ \alpha_2 =\beta_2=0$. The same contradiction is reached when trying to solve the system (\ref{eq34}). So, $ \det( \mathcal{S}) \ne 0$ in all cases. Consequently, the analogue boundary value problem (\ref{eq23})-(\ref{eq26}) in the Laplace domain, or equivalently the problem (\ref{eq01})-(\ref{eq04}) in the time domain, admits a unique solution whenever $ \alpha_1^2+\beta_1^2 \ne 0$ and $ \alpha_2^2+\beta_2^2 \ne 0$. Now, the exact solution in the Laplace domain of the problem (\ref{eq01})-(\ref{eq04}) is $ U(x,p) $ given by equation (\ref{eq27}), with the functions $ U(l_1, p) $, $ U(l_2, p) $, $ U_x(l_1, p) $ and $ U_x(l_2, p) $ expressed by using determinants as:
\begin{equation}
	U(l_1, p)=\frac{\left| 
		\begin{array}{cccc}
			\displaystyle G_1(p)	& \displaystyle\beta_1 & \displaystyle 0 &  \displaystyle 0\\
			\displaystyle G_2(p)	& 0 &\displaystyle \alpha_2 &  \displaystyle\beta_2\\
			\displaystyle R(l_1,p)	& \displaystyle \frac{a}{2 \sqrt{b+p}} &  \displaystyle\frac{-\chi(l_2-l_1,p)}{2} & \displaystyle \frac{-a\chi(l_2-l_1,p)}{2 \sqrt{b+p}} \\
			\displaystyle R(l_2,p)& \displaystyle\frac{a\chi(l_2-l_1,p)}{2 \sqrt{b+p}} & \displaystyle\frac{1}{2} & \displaystyle \frac{-a}{2 \sqrt{b+p}}
		\end{array}
		\right|}{\mathrm{det}( \mathcal{S})};  
	\label{eq36}
\end{equation}
\begin{equation}
	U_x(l_1, p)=\frac{\left| 
		\begin{array}{cccc}
			\displaystyle \alpha_1	& \displaystyle G_1(p) & \displaystyle 0 &  \displaystyle 0\\
			\displaystyle 0	& G_2(p) &\displaystyle \alpha_2 &  \displaystyle\beta_2\\
			\displaystyle\frac{1}{2}	& \displaystyle R(l_1,p) &  \displaystyle\frac{-\chi(l_2-l_1,p)}{2} & \displaystyle \frac{-a\chi(l_2-l_1,p)}{2 \sqrt{b+p}} \\
			\displaystyle\frac{-\chi(l_2-l_1,p)}{2}	& \displaystyle R(l_2,p) & \displaystyle\frac{1}{2} & \displaystyle \frac{-a}{2 \sqrt{b+p}}
		\end{array}
		\right|}{\mathrm{det}( \mathcal{S})};  
	\label{eq37}
\end{equation}
\begin{equation}
	U(l_2, p)=\frac{\left| 
		\begin{array}{cccc}
			\displaystyle \alpha_1	& \displaystyle\beta_1 & \displaystyle G_1(p) &  \displaystyle 0\\
			\displaystyle 0	& 0 &\displaystyle G_2(p) &  \displaystyle\beta_2\\
			\displaystyle\frac{1}{2}	& \displaystyle \frac{a}{2 \sqrt{b+p}} &  \displaystyle R(l_1,p) & \displaystyle \frac{-a\chi(l_2-l_1,p)}{2 \sqrt{b+p}} \\
			\displaystyle\frac{-\chi(l_2-l_1,p)}{2}	& \displaystyle\frac{a\chi(l_2-l_1,p)}{2 \sqrt{b+p}} & \displaystyle R(l_2,p) & \displaystyle \frac{-a}{2 \sqrt{b+p}}
		\end{array}
		\right|}{\mathrm{det}( \mathcal{S})};  
	\label{eq38}
\end{equation}
and 
\begin{equation}
	U_x(l_x, p)=\frac{\left| 
		\begin{array}{cccc}
			\displaystyle \alpha_1	& \displaystyle\beta_1 & \displaystyle 0 &  \displaystyle G_1(p)\\
			\displaystyle 0	& 0 &\displaystyle \alpha_2 &  \displaystyle G_2(p)\\
			\displaystyle\frac{1}{2}	& \displaystyle \frac{a}{2 \sqrt{b+p}} &  \displaystyle\frac{-\chi(l_2-l_1,p)}{2} & \displaystyle R(l_1,p) \\
			\displaystyle\frac{-\chi(l_2-l_1,p)}{2}	& \displaystyle\frac{a\chi(l_2-l_1,p)}{2 \sqrt{b+p}} & \displaystyle\frac{1}{2} & \displaystyle R(l_2,p)
		\end{array}
		\right|}{\mathrm{det}( \mathcal{S})}.  
	\label{eq39}
\end{equation}

In brief, the exact solution in the  Laplace domain of the boundary value problem (\ref{eq01})-(\ref{eq04}) is expressed in a unified way by the function $ U(x,p) $ given by equation (\ref{eq27}) together with those given in equations (\ref{eq36})-(\ref{eq39}). Exact series solutions in the time domain of such linear boundary value problems are well established by using the Fourier decomposition method as in \cite{An-series} for example. Thus, the exact Laplace transform of those series solutions performed via the Sturm-Liouville theory can be recovered by the expression of the function $ U(x,p) $.       

\section{Analytical approximations in short time limits}
\label{sec:5}
In this section, our interest is to find approximate analytical solutions to the problem (\ref{eq01})-(\ref{eq04}) at the earliest times of the process. We assume a subdivision of the time interval $ [0, T] $, such that the magnitude of the dimensionless time step $ \Delta t$ is sufficiently small, for example $ \Delta t \le 10^{-2}$, as in many schemes in numerical analysis (see \cite{Su-Ji-Ji-numerical} and \cite{Iz-Yu-hybrid} among others). The goal of this section is to derive from the exact operational solution (\ref{eq27}), approximate analytical solutions valid during the first time step of the reaction–diffusion process, namely for $t\in  [0, \Delta t]$. 

The limiting case of a short time ($ \Delta t$ tending to 0) corresponds to a very large value of the Laplace domain variable ($ p $ tending to $+\infty$). In order to minimize calculations, some prior simplifications can be done on the determinant of the system $ \mathrm{det}( \mathcal{S})$, before deriving the asymptotic expansions of the solution and the related truncated expansions in the time domain. For $ \eta>0 $, the inverse Laplace transform of $ \chi(\eta,p) $ is written as:
\begin{equation}
	{\mathcal L}^{-1}\{\chi(\eta,p)\}={\mathcal L}^{-1}\left\lbrace \exp\left(\frac{-\eta\sqrt{b+p}}{a} \right)\right\rbrace =\frac{1}{2}\frac{\eta \exp\left( -bt-\frac{\eta^2}{4a^2 t}\right) }{a\sqrt{\pi}t^{3/2}} =\Gamma(t). 
	\label{eq40}
\end{equation}
When $ p\rightarrow \infty $, $ \chi(\eta,p) $ is negligible compared to $ {1}/{p^n} $, for all positive integers $ n\geq 1 $, whereas when $ t\rightarrow 0 $, the inverse Laplace transform $ \Gamma(t) $ is negligible compared to $ t^n $. Indeed, $\displaystyle\lim_{p\rightarrow \infty} p^n\chi(\eta,p)=\displaystyle\lim_{t\rightarrow 0} \Gamma(t)/t^n=0$, implying $ \chi(\eta,p)=o({1}/{p^n}) $ at $ p=\infty $ and $  \Gamma(t)=o(t^n) $ at $ t=0 $, where the Little-o is the asymptotic notation. Again, these relations are valid for arbitrary order $n\geq 1$. When $ p\rightarrow \infty $, the term $ \chi(l_2-l_1,p) $ and its asymptotic expansions are negligible in the Laplace domain as well as their inverses in the time domain when $ t\rightarrow 0 $. Therefore, when $ p $ is sufficiently large, the asymptotic solution can be obtained using the following determinant for the system formed by equations (\ref{eq25}), (\ref{eq26}),  (\ref{eq28}) and (\ref{eq29}): 
\begin{equation}
	\mathrm{det}( \mathcal{S}^a)=\left| 
	\begin{array}{cccc}
		\displaystyle \alpha_1	& \displaystyle\beta_1 & \displaystyle 0 &  \displaystyle 0\\
		\displaystyle 0	& 0 &\displaystyle \alpha_2 &  \displaystyle\beta_2\\
		\displaystyle\frac{1}{2}	& \displaystyle \frac{a}{2 \sqrt{b+p}} &  \displaystyle 0 & \displaystyle 0 \\
		\displaystyle 0	& \displaystyle 0 & \displaystyle\frac{1}{2} & \displaystyle \frac{-a}{2 \sqrt{b+p}}
	\end{array}
	\right|.  
	\label{eq41}
\end{equation}

Based on the same assumption that the function $ \chi(\eta,p) $ is negligible for all $ \eta>0 $ whenever $ p $ is sufficiently large, an asymptotic expansion of the exact solution $ U(x,p) $, valid at any order $ n\geq 1 $, can be deduced for $l_1<x< l_2 $, from equation (\ref{eq27}) as follows: 
\begin{equation}
	U^a(x,p)=R(x,p)+o\left( \frac{1}{p^{n}} \right).	  
	\label{eq42}
\end{equation}
Equation (\ref{eq42}) corresponds for $0<t< \Delta t$ to the following truncated expansion of the solution $ 	u(x,t) $ in the time domain:
\begin{equation}
	u^a(x,t)=r(x,t)+o\left( (\Delta t)^{n} \right);	  
	\label{eq43}
\end{equation}
where $ r(x,t) $ is given by equation (\ref{eq20}) and $ o\left( (\Delta t)^{n} \right) $  is an upper bound of the truncation error $ o\left(t^{n} \right) $. But, at $ x=l_1 $ and $ x=l_2 $, some refined asymptotic expansions in the Laplace domain $ U^a(l_1, p) $, $ U^a(l_2, p) $, $ U^a_x(l_1, p) $ and $ U^a_x(l_2, p) $ can be obtained for the boundary-related functions $ U(l_1, p) $, $ U(l_2, p) $, $ U_x(l_1, p) $ and $ U_x(l_2, p) $. We distinguish four different cases with respect to the values of the two coefficients $ \beta_1$ and $\beta_2$. For each case, the corresponding analytical approximations in the time domain, namely  $ u^a(l_1, t) $, $ u^a(l_2, t) $, $ u^a_x(l_1, t) $ and $ u^a_x(l_2, t) $, are also given for $t\in [0, \Delta t]$. Although this can be improved, the asymptotic expansions at second order are retained here, while the corresponding analytical approximations in the time domain are given at the first order.

	\vspace{2cm}
\begin{itemize}
	\item \underline{Case $ \beta_1 \beta_2\ne 0$}
	\begin{enumerate}
		\item Solving the system by using the reduced determinant $ \mathrm{det}( \mathcal{S}^a) $ expressed above in formula (\ref{eq41}), leads to:
		$$ \displaystyle U(l_1, p)= -\frac{a}{\beta_1 \sqrt{b + p}- a \alpha_1  }\,G_1(p)+ \frac{2 \sqrt{b + p}  \beta_1}{\beta_1 \sqrt{b + p}- a \alpha_1 }\,R(l_1, p),$$ and if the Big-O denotes the asymptotic notation, it can be written:
		$$
		\displaystyle-\frac{a}{\beta_1 \sqrt{b + p}- a \alpha_1  }=U^{1l_1}(p)+O \left(  \frac{1}{p^{2}} \right);
		$$
		and
		$$
		\displaystyle\frac{2 \sqrt{b + p} \beta_1}{\beta_1 \sqrt{b + p}- a \alpha_1 }=2+U^{2l_1}(p)+O \left(  \frac{1}{p^{2}} \right);
		$$
		where
		$$
		U^{1l_1}(p)=\displaystyle -\frac {a }{\beta_1}\sqrt{\frac{1}{p}}-{\frac {{a}^{2}\alpha_1
			}{{\beta_1}^{2}p}}-{\frac {1}{\beta_
				1} \left( -\frac{1}{2}\,ab+{\frac {{a}^{3}{\alpha_1}^{2}}{{\beta_1
					}^{2}}} \right)\left( \frac{1}{p} \right) ^{3/2} };
		$$
		and
		$$
		U^{2l_1}(p)=\displaystyle 2\,{\frac {a\alpha_1 }{\beta_1}}\sqrt{\frac{1}{p}}+2\,\frac {{
				a}^{2}{\alpha_1}^{2}}{{\beta_1}^{2}p}+\,{\frac {2}{\beta_1} \left( -\frac{1}{2}\,a\alpha_1 b+{\frac {
					{a}^{3}{\alpha_1}^{3}}{{\beta_1}^{2}}} \right) \left( \frac{1}{p} \right) ^{3/2} }.
		$$
		Now, due to the properties of the LIT, $ G_1(p) $ and $R(l_1, p)  $ are bounded functions in the Laplace domain for $ p>0 $, and an asymptotic expansion of the solution $ U(l_1, p) $ can be written as:
		$$
		\displaystyle U^a(l_1, p)=U^{1l_1}(p)G_1(p)+U^{2l_1}(p)R(l_1, p)+2R(l_1, p)+O \left(  \frac{1}{p^{2}} \right).
		$$ 
		The corresponding analytical approximation in the time domain during a short time step $t\in [0, \Delta t]$ reads 
		$$u^a(l_1, t)= u^{1l_1}(t)\ast g_1(t)+u^{2l_1}(t)\ast r(l_1, t)+2r(l_1, t)+O \left( \Delta t \right), $$ where $ O \left( \Delta t \right) $ instead of $O\left( t \right) $ is an upper bound of the truncation error of the approximation, $ \ast $ denotes the convolution product, $ r(l_1, t) $ is calculated through equation (\ref{eq20}), $ g_1(t) $ is the known function related to the boundary $ x=l_1 $, 
		$$ u^{1l_1}(t)={\mathcal L}^{-1}\{U^{1l_1}(p)\}=-{\frac {{a}^{2}\alpha_1}{{\beta_1}^{2}}}-{\frac {a}{ \sqrt{
					\pi \,t}\beta_1}}+{\frac { \left( -2\,{a}^{2}{\alpha_1}^{2}+b{
					\beta_1}^{2} \right) a}{{\beta_1}^{3}} \sqrt{{\frac {t}{\pi }}
		}}; 
		$$
		and
		$$ u^{2l_1}(t)={\mathcal L}^{-1}\{U^{2l_1}(p)\}=2\,{\frac {a\alpha_1}{ \sqrt{\pi \,t}\beta_1}}+2\,{\frac {{a}^
				{2}{\alpha_1}^{2}}{{\beta_1}^{2}}}+2\,{\frac {a\alpha_1
				\left( 2\,{a}^{2}{\alpha_1}^{2}-b{\beta_1}^{2} \right) }{{
					\beta_1}^{3}} \sqrt{{\frac {t}{\pi }}}}. 
		$$
		\item For $ U(l_2,p) $ and  $ u(l_2,t) $, similar calculations as above give:
		$$
		U(l_2,p)={\frac {G_2 \left( p \right) a}{\beta_2 \sqrt{b+p}+\alpha_{{2}
				}a}}+2\,{\frac {R \left( l_2,p \right)  \sqrt{b+p}\beta_2}{
				\beta_2 \sqrt{b+p}+\alpha_2 a}};
		$$
		and
		$$
		\displaystyle U^a(l_2, p)=U^{1l_2}(p)G_2(p)+U^{2l_2}(p)R(l_2, p)+2R(l_2, p)+O \left(  \frac{1}{p^{2}} \right);
		$$ 
		where
		$$
		U^{1l_2}(p)=\displaystyle{\frac {a }{\beta_2}}\sqrt{\frac{1}{p}}-{\frac {{a}^{2}\alpha_2}{
				{\beta_2}^{2}p}}+{\frac { 1}{\beta_{{
						2}}} \left( -\frac{1}{2}\,ab+{\frac {{a}^{3}{\alpha_2}^{2}}{{\beta_2}^
					{2}}} \right)\left( \frac{1}{p} \right) ^{3/2} };
		$$
		and
		$$
		U^{2l_2}(p)=\displaystyle -2\,{\frac {\alpha_2a }{\beta_2}}\sqrt{\frac{1}{p}}+2\,{\frac {{
					a}^{2}{\alpha_2}^{2}}{{\beta_2}^{2}p}}+\,{\frac {2}{\beta_2} \left( \frac{1}{2}\,\alpha_2ab-{\frac {{
						\alpha_2}^{3}{a}^{3}}{{\beta_2}^{2}}} \right) \left(\frac{1}{p} \right) ^{3/2} }.
		$$
		The corresponding analytical approximation in the time domain during a short time step $t\in [0, \Delta t]$ is: 
		$$u^a(l_2, t)= u^{1l_2}(t)\ast g_2(t)+u^{2l_2}(t)\ast r(l_2, t)+2r(l_2, t)+O \left( \Delta t \right), $$
		where
		$$ u^{1l_2}(t)={\mathcal L}^{-1}\{U^{1l_2}(p)\}=-{\frac {{a}^{2}\alpha_2}{{\beta_2}^{2}}}+{\frac {a}{ \sqrt{
					\pi \,t}\beta_2}}+{\frac {a \left( 2\,{a}^{2}{\alpha_2}^{2}-b{
					\beta_2}^{2} \right) }{{\beta_2}^{3}} \sqrt{{\frac {t}{\pi }}}}; 
		$$
		and
		$$ u^{2l_2}(t)={\mathcal L}^{-1}\{U^{2l_2}(p)\}=-2\,{\frac {\alpha_2a}{ \sqrt{\pi \,t}\beta_2}}+2\,{\frac {{a}
				^{2}{\alpha_2}^{2}}{{\beta_2}^{2}}}+2\,{\frac { \left( -2\,{a}
				^{2}{\alpha_2}^{2}+b{\beta_2}^{2} \right) \alpha_2a}{{
					\beta_2}^{3}} \sqrt{{\frac {t}{\pi }}}}. 
		$$
		
		\item Likewise, for $ U_x(l_1,p) $ and  $ u_x(l_1,t) $, the results are:
		$$
		U_x(l_1,p)={\frac { \sqrt{b+p}\,G_1 \left( p \right) }{\beta_1 \sqrt{b+p}-a
				\alpha_1}}-2\,{\frac {\alpha_1 \sqrt{b+p}\,R \left( l_1,p
				\right) }{\beta_1 \sqrt{b+p}-a\alpha_1}};
		$$
		and
		$$
		\displaystyle U_x^a(l_1, p)={\beta_1}^{-1}\,G_1(p)+U_x^{1l_1}(p)G_1(p)+U_x^{2l_1}(p)R(l_1, p)-2\,{\frac {\alpha_1}{\beta_1}}\,R(l_1, p)+O \left(  \frac{1}{p^{2}} \right);
		$$ 
		where
		$$
		U_x^{1l_1}(p)=\displaystyle{\frac {a\alpha_1}{{\beta_1
				}^{2}}} \sqrt{\frac{1}{p}}+{\frac {{a}^{2}{\alpha_1}^{2}}{{\beta_1}^{3}p}}+{
			\frac { 1}{\beta_1} \left( -\frac{1}{2}\,{
				\frac {a\alpha_1b}{\beta_1}}+{\frac {{a}^{3}{\alpha_1}^{3}
				}{{\beta_1}^{3}}} \right)\left( \frac{1}{p} \right) ^{3/2} };
		$$
		and
		$$
		U_x^{2l_1}(p)=\displaystyle -2\,{\frac {{\alpha_1}^{2}a
			}{{\beta_1}^{2}}}\sqrt{\frac{1}{p}}-2\,{\frac {{\alpha_1}^{3}{a}
				^{2}}{{\beta_1}^{3}p}}-\,{\frac {2 
			}{\beta_1} \left( -\frac{1}{2}\,{\frac {{\alpha_1}^{2}ab}{\beta_1}
			}+{\frac {{\alpha_1}^{4}{a}^{3}}{{\beta_1}^{3}}} \right)\left( \frac{1}{p} \right) ^{3/2} }.
		$$
		The related analytical approximation in the time domain during a short time step $t\in [0, \Delta t]$ is: 
		$$u_x^a(l_1, t)= {\beta_1}^{-1}\,g_1(t)+u_x^{1l_1}(t)\ast g_1(t)+u_x^{2l_1}(t)\ast r(l_1, t)-2\,{\frac {\alpha_1}{\beta_1}}\,r(l_1, t)+O \left( \Delta t \right), $$
		where
		$$ u_x^{1l_1}(t)={\mathcal L}^{-1}\{U_x^{1l_1}(p)\}={\frac {a\alpha_1}{ \sqrt{\pi \,t}{\beta_1}^{2}}}+{\frac {{a}^
				{2}{\alpha_1}^{2}}{{\beta_1}^{3}}}+{\frac {a\alpha_1
				\left( 2\,{a}^{2}{\alpha_1}^{2}-b{\beta_1}^{2} \right) }{{
					\beta_1}^{4}} \sqrt{{\frac {t}{\pi }}}}; 
		$$
		and
		$$ u_x^{2l_1}(t)={\mathcal L}^{-1}\{U_x^{2l_1}(p)\}=-2\,{\frac {{\alpha_1}^{2}a}{ \sqrt{\pi \,t}{\beta_1}^{2}}}-2
		\,{\frac {{\alpha_1}^{3}{a}^{2}}{{\beta_1}^{3}}}+2\,{\frac {
				\left( -2\,{a}^{2}{\alpha_1}^{2}+b{\beta_1}^{2} \right) {
					\alpha_1}^{2}a}{{\beta_1}^{4}} \sqrt{{\frac {t}{\pi }}}}. 
		$$
		
		\item Finally for $ U_x(l_2,p) $ and  $ u_x(l_2,t) $, calculations give:
		$$
		U_x(l_2,p)={\frac { \sqrt{b+p}\,G_2 \left( p \right) }{\beta_2 \sqrt{b+p}+
				\alpha_2a}}-2\,{\frac {\alpha_2 \sqrt{b+p}\,R \left( l_2,p
				\right) }{\beta_2 \sqrt{b+p}+\alpha_2a}};
		$$
		and
		$$
		\displaystyle U_x^a(l_2, p)={\beta_2}^{-1}\,G_2(p)+U_x^{1l_2}(p)G_2(p)+U_x^{2l_2}(p)R(l_2, p)-2\,{\frac {\alpha_2}{\beta_2}}\,R(l_2, p)+O \left(  \frac{1}{p^{2}} \right);
		$$ 
		where
		$$
		U_x^{1l_2}(p)=\displaystyle-{\frac {\alpha_2a }{{\beta_2
				}^{2}}}\sqrt{\frac{1}{p}}+{\frac {{a}^{2}{\alpha_2}^{2}}{{\beta_2}^{3}p}}+{
			\frac { 1}{\beta_2} \left( \frac{1}{2}\,{
				\frac {\alpha_2ab}{\beta_2}}-{\frac {{\alpha_2}^{3}{a}^{3}
				}{{\beta_2}^{3}}} \right)\left( \frac{1}{p} \right) ^{3/2} };
		$$
		and
		$$
		U_x^{2l_2}(p)=\displaystyle 2\,{\frac {{\alpha_2}^{2}a
			}{{\beta_2}^{2}}}\sqrt{\frac{1}{p}}-2\,{\frac {{\alpha_2}^{3}{a}
				^{2}}{{\beta_2}^{3}p}}-\,{\frac {2 }{\beta_2} \left( \frac{1}{2}\,{\frac {{\alpha_2}^{2}ab}{\beta_2}}
			-{\frac {{\alpha_2}^{4}{a}^{3}}{{\beta_2}^{3}}} \right)\left( \frac{1}{p} \right) ^{3/2} }.
		$$
		The corresponding analytical truncation in the time domain during a short time step $t\in [0, \Delta t]$ is written as: 
		$$u_x^a(l_2, t)= {\beta_2}^{-1}\,g_2(t)+u_x^{1l_2}(t)\ast g_2(t)+u_x^{2l_2}(t)\ast r(l_2, t)-2\,{\frac {\alpha_2}{\beta_2}}\,r(l_2, t)+O \left( \Delta t \right), $$
		where
		$$ u_x^{1l_2}(t)={\mathcal L}^{-1}\{U_x^{1l_2}(p)\}=-{\frac {\alpha_2a}{ \sqrt{\pi \,t}{\beta_2}^{2}}}+{\frac {{a}
				^{2}{\alpha_2}^{2}}{{\beta_2}^{3}}}+{\frac { \left( -2\,{a}^{2
				}{\alpha_2}^{2}+b{\beta_2}^{2} \right) \alpha_2a}{{\beta_{
						{2}}}^{4}} \sqrt{{\frac {t}{\pi }}}}; 
		$$
		and
		$$ u_x^{2l_2}(t)={\mathcal L}^{-1}\{U_x^{2l_2}(p)\}=2\,{\frac {{\alpha_2}^{2}a}{ \sqrt{\pi \,t}{\beta_2}^{2}}}-2\,
		{\frac {{\alpha_2}^{3}{a}^{2}}{{\beta_2}^{3}}}+2\,{\frac {{
					\alpha_2}^{2}a \left( 2\,{a}^{2}{\alpha_2}^{2}-b{\beta_2}^
				{2} \right) }{{\beta_2}^{4}} \sqrt{{\frac {t}{\pi }}}}. 
		$$
	\end{enumerate}

	\vspace{2cm}

	\item \underline{Case $ \beta_1\ne 0$ and $ \beta_2= 0$}
	\begin{enumerate}
		\item As in the first case, the system is solved using the reduced determinant $ \mathrm{det}( \mathcal{S}^a) $. This leads to:
		$$ \displaystyle U(l_1, p)= -{\frac {G_1 \left( p \right) a}{\beta_1 \sqrt{b+p}-a\alpha_{{
						1}}}}+2\,{\frac { \sqrt{b+p}\,R \left( l_1,p \right) \beta_1}{
				\beta_1 \sqrt{b+p}-a\alpha_1}},$$ 
		and $ U(l_1,p) $ and  its asymptotic expansion $ U^a(l_1,p) $, as well as $ u(l_1,t)$ and its truncation expansion $ u^a(l_1,t)$ are the same as in the case $ \beta_1\beta_2 \ne 0 $, item (1).

		\item Concerning $ U(l_2,p) $ and  $ u(l_2,t) $, the results are reduced to:
		$$
		U(l_2,p)={\frac {G_2 \left( p \right) }{\alpha_2}};
		$$
		and the reverse in the time domain during a short time step $t\in [0, \Delta t]$ is: 
		$$u(l_2, t)= \frac{g_2(t)}{\alpha_2}.$$

		\item The functions $ U_x(l_1,p) $,  $ u_x(l_1,t)$ and their asymptotic and truncation expansions $ U^a_x(l_1,p) $,  $ u^a_x(l_1,t)$ are the same as in the case $ \beta_1\beta_2 \ne 0 $, item (3), since
		$$
		U_x(l_1,p)={\frac { \sqrt{b+p}\,G_1 \left( p \right) }{\beta_1 \sqrt{b+p}-a
				\alpha_1}}-2\,{\frac {\alpha_1 \sqrt{b+p}\,R \left( l_1,p
				\right) }{\beta_1 \sqrt{b+p}-a\alpha_1}}.
		$$

		\item About $ U_x(l_2,p) $ and  $ u_x(l_2,t) $, one has:
		$$
		U_x(l_2,p)={\frac { \sqrt{b+p}\,G_2 \left( p \right) }{\alpha_2 a}}-2\,{
			\frac { \sqrt{b+p}\,R \left( l_2,p \right) }{a}};
		$$
		and the calculations will use here the assumptions that $ r(l_2,t) $ and $ g_2(t) $ are once differentiable relatively to the time variable $ t $ on $ t\in [0, \Delta t] $, and their respective derivative $ r'(l_2,t) $ and $ g'_2(t) $ verify the relations: $r(l_2,t)-r(l_2,0)=\displaystyle\int_{0}^{t}\,r'(l_2,\tau)d\tau
		$, and $g_2(t)-g_2(0)=\displaystyle\int_{0}^{t}\,g'_2(\tau)d\tau.
		$ The derivatives are assumed to be of exponential order so that if $ R^q(l_2,p) $ and $ G^q_2(p) $ are their Laplace transforms, we have according to the properties: 
		$$R(l_2,p)=\displaystyle{\mathcal L}\left\lbrace \int_{0}^{t}\,r'(l_2,\tau)d\tau\right\rbrace +\frac{r(l_2,0)}{p} =\frac{R^q(l_2,p)}{p}+\frac{r(l_2,0)}{p} 
		$$
		and
		$$G_2(p)=\displaystyle{\mathcal L}\left\lbrace \int_{0}^{t}\,g'_2(\tau)d\tau\right\rbrace+\frac{g_2(0)}{p}  =\frac{G^q_2(p)}{p}+\frac{g_2(0)}{p} .
		$$
		Then
		$$
		U_x(l_2,p)={\frac { \sqrt{b+p}\,\left( G^q_2(p)+g_2(0) \right) }{\alpha_2a \,p}}-2\,{
			\frac { \sqrt{b+p}\,\left( R^q (l_2,p )+r(l_2,0) \right) }{a\,p}};
		$$
		and
		$$
		\displaystyle U_x^a(l_2, p)=U_x^{l_2}(p)\left( \frac{ G^q_2(p)+g_2(0) }{a\alpha_2}-2\frac{ R^q (l_2,p )+r(l_2,0) }{a}\right) +O \left(  \frac{1}{p^{5/2}} \right);
		$$ 
		where
		$$
		U_x^{l_2}(p)=\displaystyle \sqrt{\frac{1}{p}}+\,{\frac {b }{2}}\left( \frac{1}{p} \right) ^{3/2}. 
		$$
		In this sub-case, a truncation expansion in the time domain during a short time step $t\in [0, \Delta t]$ is: 
		$$u_x^a(l_2, t)= \frac{1}{a}\left(  \frac{g_2(0)}{\alpha_2}- 2r(l_2, 0)\right)u_x^{l_2}(t) +\frac{1}{a}u_x^{l_2}(t)\ast\left(  \frac{g'_2(t)}{\alpha_2}- 2r'(l_2, t)\right) +O \left( \Delta t \right), $$
		where
		$$ u_x^{l_2}(t)={\frac {bt+1}{ \sqrt{\pi \,t}}}. 
		$$
	\end{enumerate}
	
	\vspace{2cm}
	
	\item \underline{Case $ \beta_1= 0$ and $ \beta_2\ne 0$}
	\begin{enumerate}
		
		\item The function $ U(l_1,p) $ and  its reverse $ u(l_1,t) $ in the time domain ($t\in [0, \Delta t]$) are respectively:
		$$
		U(l_1,p)={\frac {G_1 \left( p \right) }{\alpha_1}};
		$$
		and
		$$u(l_1, t)= \frac{g_1(t)}{\alpha_1}.$$
		
		\item For $ U(l_2,p) $ and  $ u(l_2,t) $, the results are identical to those obtained in the case $ \beta_1 \beta_2\ne 0 $, item (2) since:
		$$
		U(l_2,p)={\frac {G_2\left( p \right) a}{\beta_2 \sqrt{b+p}+\alpha_{{2}
				}a}}+2\,{\frac {R \left( l_2,p \right)  \sqrt{b+p}\beta_2}{
				\beta_2 \sqrt{b+p}+\alpha_2 a}};
		$$
		
		\item For $ U_x(l_1,p) $ and  $ u_x(l_1,t) $, calculations give:
		$$
		U_x(l_1,p)=-{\frac { \sqrt{b+p}\,G_1 \left( p \right) }{\alpha_1a}}+2\,{
			\frac { \sqrt{b+p}\,R \left( l_1,p \right) }{a}};
		$$
		and the assumptions are that $ r(l_1,t) $ and $ g_1(t) $ are once differentiable in rapport to the time variable $ t $ on $ t\in [0, \Delta t] $. The respective derivatives $ r'(l_1,t) $ and $ g'_1(t) $ can be related to the original functions by the relations: $r(l_1,t)-r(l_1,0)=\displaystyle\int_{0}^{t}\,r'(l_1,\tau)d\tau
		$, and $g_1(t)-g_1(0)=\displaystyle\int_{0}^{t}\,g'_1(\tau)d\tau.$ Since the derivatives are assumed to be of exponential order, their Laplace transforms $ R^q(l_1,p) $ and $ G^q_1(p) $ verify the properties: 
		$$R(l_1,p)=\displaystyle{\mathcal L}\left\lbrace \int_{0}^{t}\,r'(l_1,\tau)d\tau\right\rbrace +\frac{r(l_1,0)}{p} =\frac{R^q(l_1,p)}{p}+\frac{r(l_1,0)}{p} 
		$$
		and
		$$G_1(p)=\displaystyle{\mathcal L}\left\lbrace \int_{0}^{t}\,g'_1(\tau)d\tau\right\rbrace+\frac{g_1(0)}{p}  =\frac{G^q_1(p)}{p}+\frac{g_1(0)}{p} .
		$$
		Then
		$$
		U_x(l_1,p)=-{\frac { \sqrt{b+p}\,\left( G^q_1(p)+g_1(0) \right) }{\alpha_1a \,p}}+2\,{
			\frac { \sqrt{b+p}\,\left( R^q (l_1,p )+r(l_1,0) \right) }{a\,p}};
		$$
		and
		$$
		\displaystyle U_x^a(l_1, p)=U_x^{l_1}(p)\left( -\frac{ G^q_1(p)+g_1(0) }{a\alpha_1}+2\frac{ R^q (l_1,p )+r(l_1,0) }{a}\right) +O \left(  \frac{1}{p^{5/2}} \right);
		$$ 
		where
		$$
		U_x^{l_1}(p)=\displaystyle \sqrt{\frac{1}{p}}+\,{\frac {b }{2}}\left( \frac{1}{p} \right) ^{3/2}. 
		$$
		A truncation expansion in the time domain during a short time step $t\in [0, \Delta t]$ corresponds to: 
		$$u_x^a(l_1, t)= \frac{1}{a}\left(-  \frac{g_1(0)}{\alpha_1}+ 2r(l_1, 0)\right)u_x^{l_1}(t) +\frac{1}{a}u_x^{l_1}(t)\ast\left(-  \frac{g'_1(t)}{\alpha_1}+ 2r'(l_1, t)\right) +O \left( \Delta t \right), $$
		where
		$$ u_x^{l_1}(t)={\frac {bt+1}{ \sqrt{\pi \,t}}}. 
		$$
		
		\item Since:
		$$
		U_x(l_2,p)={\frac { \sqrt{b+p}\,G_2 \left( p \right) }{\beta_2 \sqrt{b+p}+
				\alpha_2a}}-2\,{\frac {\alpha_2 \sqrt{b+p}\,R \left( l_2,p
				\right) }{\beta_2 \sqrt{b+p}+\alpha_2a}};
		$$
		the results for $ U_x(l_2,p) $ and  $ u_x(l_2,t) $ in this sub-case are identical to those obtained in the case $ \beta_1 \beta_2\ne 0 $, item (4). 
		
	\end{enumerate}
	
	\vspace{2cm}
	
	\item \underline{Case $ \beta_1= 0$ and $ \beta_2= 0$}
	\begin{enumerate}
		\item The expression of $ U(l_1,p) $ and of its reverse $ u(l_1,t) $ are identical to those obtained in the case $ \beta_1= 0$ and $ \beta_2\ne 0$, item (1).
		
		\item The expression of $ U(l_2,p) $ and of its reverse $ u(l_2,t) $ are the same as those obtained in the case $ \beta_1\ne 0$ and $ \beta_2= 0$, item (2).
		
		\item The expansions of $ U_x(l_1,p) $ and of its reverse $ u_x(l_1,t) $ are equal to those obtained in the case $ \beta_1= 0$ and $ \beta_2\ne 0$, item (3).
		
		\item The expansions of $ U_x(l_2,p) $ and of its reverse $ u_x(l_2,t) $ are identical those obtained in the case $ \beta_1\ne 0$ and $ \beta_2= 0$, item (4).
	\end{enumerate}

\end{itemize}

\section{Discussion of the results and applications}
\label{sec:6} 
First, the results mentioned in this document can be compared and discussed on a specific example. Some powerful PDE toolbox functions exist in Matlab software for example, and scripts based on Gauss-Seidel and finite difference methods are available online. Therefore, the curves obtained from the Fourier decomposition method, when using partial sums of the infinite series solution, and those obtained from the method of approximate analytical solutions in short time limits developed in this paper, can both be compared to numerical solution curves, when using toolbox functions. For this purpose, the example titled Example 6.1. is considered from \cite{Hen-partial}. It consists to solve a problem where the initial condition matches with the boundary conditions. The equation is stated as:
\begin{equation}
	\frac{\partial u}{\partial t}-a^2\frac{\partial^{2}u}{\partial x^{2}} =0, \quad 0 < x < l,\quad  t>0,
	\label{eq44}
\end{equation}    
subject to consistent initial and boundary conditions
\begin{equation}
	u(x,0)=\varphi(x), \quad \quad  u(0,t)=u(0,t)=0, 
	\label{eq45}
\end{equation}
where
\begin{equation}
	\varphi(x)=\left\lbrace 
	\begin{array}{cccc}
		\displaystyle\frac{x}{l}u_0\quad \mathrm{for}\quad 0 \le x \le \frac{l}{2}\\\\
		\displaystyle\frac{l-x}{l}u_0\quad{for} \quad  \frac{l}{2}< x \le l
	\end{array}
	\right. 
	\label{eq46}
\end{equation}
with $ u_0 $ being a constant. The exact solution obtained by the Fourier decomposition method can be reported as the following infinite series:
\begin{equation}
	\displaystyle u(x,t)=\frac{4 u_0}{\pi^2}\sum_{k=1}^{+\infty}\frac{(-1)^{k+1}}{(2k-1)^2}\exp\left( -\frac{a^2(2k-1)^2\pi^2}{l^2}t\right) \sin \frac{(2k-1)\pi x}{l}.
	\label{eq47}
\end{equation}
Assigning the corresponding parameter values to the general problem (\ref{eq01})-(\ref{eq04}), that is, $ l_1 =0,\  l_2=l, \ b=0, \ f(x,t)=0, \ \beta_1=\beta_2=0, \ \alpha_1=\alpha_2 =1, \ g_1(t)=g_2(t)=0,$ the p-domain solution (\ref{eq27}) is reduced to:
\begin{equation}
	\begin{array}{ll} 
		U(x,p)=	\displaystyle \frac{a}{2 \sqrt{p}}\left[U_x(l,p)\exp\left(\frac{-(l-x)\sqrt{p}}{a} \right)-U_x(0,p)\exp\left(\frac{-x\sqrt{p}}{a} \right)   \right]\\ 
		+\displaystyle{\mathcal L}\left\lbrace \frac{1}{2 a\sqrt{\pi t}}
		\int_{0}^{l}\varphi(\xi)\left[\exp\left( -\frac{(\xi-x)^2}{4 a^2 t}\right) \right] d\xi\right\rbrace.
	\end{array} 
	\label{eq48}
\end{equation}
Then, for the case when $l = 10$, $u_0 = 5$ and $a^2 = 0.25$ as in the Example 6.1. from \cite{Hen-partial}, the function $ \varphi $ is rewritten as:
$$
\varphi(x)=\left\lbrace 
\begin{array}{cccc}
	\displaystyle\frac{1}{2}x\quad \mathrm{for}\quad 0 \le x \le 5\\\\
	\displaystyle 5-\frac{1}{2}x\quad{for} \quad  5 < x \le 10;
\end{array} 
\right. 
$$ 
while the reported series (\ref{eq47}) becomes:
\begin{equation}
	\displaystyle u(x,t)=\frac{20}{\pi^2}\sum_{k=1}^{+\infty}\frac{(-1)^{k+1}}{(2k-1)^2}\exp\left(-\frac{(2k-1)^2\pi^2}{400}t\right) \sin \frac{(2k-1)\pi x}{10}.
	\label{eq49}
\end{equation}
By using equations (\ref{eq37}) and (\ref{eq39}) and denoting $ \exp(x) $ by ${\rm e}^x  $,  the exact operational solution (\ref{eq48}) can be written as:
\begin{equation}
	U(x,p)=	\displaystyle \frac{1}{4 \sqrt{p}}\left[U_x(10,p)\,\exp\left(-2(10-x)\sqrt{p}\right)-U_x(0,p)\,\exp\left(-2x\sqrt{p} \right) \right] +R(x,p).
	\label{eq50}
\end{equation}
In the above expression,
\begin{equation}
	U_x(0,p)=-U_x(10,p)=\displaystyle-\frac{1}{2}\,{\frac {-{{\rm e}^{-20\, \sqrt{p}}}+2\,{{\rm e}^{-10\, \sqrt{p}}}-1}{ \left( {{\rm e}^{-20\, \sqrt{p}}}+1 \right) p}};	
	\label{eq51}	
\end{equation} 
and
$$
R(x,p)=\left\lbrace 
\begin{array}{cccc}
	\displaystyle-\frac{1}{8}\,{\frac {-4\,x \sqrt{p}+2\,{{\rm e}^{2\, \left( -5+x \right) 
					\sqrt{p}}}-{{\rm e}^{-2\,x \sqrt{p}}}-{{\rm e}^{2\, \left( -10+x
					\right)  \sqrt{p}}}}{{p}^{\frac{3}{2}}}}
	\quad \mathrm{for}\quad 0 \le x \le 5,\\\\
	\displaystyle-\frac{1}{8}\,{\frac {4\,x \sqrt{p}-40\, \sqrt{p}+2\,{{\rm e}^{-2\, \left( -5+
					x \right)  \sqrt{p}}}-{{\rm e}^{-2\,x \sqrt{p}}}-{{\rm e}^{2\, \left( 
					-10+x \right)  \sqrt{p}}}}{{p}^{\frac{3}{2}}}}
	\quad{ \rm for} \quad  5 < x \le 10.
\end{array}
\right. 
$$
The formula (\ref{eq23}) can then be checked, and the exactness of the operational solution (\ref{eq50}) is proved, that is:
$$
-a^2\frac{d^2 U}{dx^{2}}(x,p)+pU(x,p)=\varphi(x).
$$
Now, according to the results obtained in section  \ref{sec:5} relatively to the present case ($ \beta_1=\beta_2=0 $), the corresponding truncation expansion in the time domain during the short time step $t\in [0, \Delta t]$, are respectively recalled as: 
$$	u^a(x,t)=r(x,t)+o\left( (\Delta t)^{n} \right);	$$
and
$$u_x^a(l_1, t)= \frac{1}{a}\left(-  \frac{g_1(0)}{\alpha_1}+ 2r(l_1, 0)\right)u_x^{l_1}(t) +\frac{1}{a}u_x^{l_1}(t)\ast\left(-  \frac{g'_1(t)}{\alpha_1}+ 2r'(l_1, t)\right) +O \left( \Delta t \right), $$
where
$$ u_x^{l_1}(t)={\frac {bt+1}{ \sqrt{\pi \,t}}}, 
$$
and again
$$u_x^a(l_2, t)= \frac{1}{a}\left(  \frac{g_2(0)}{\alpha_2}- 2r(l_2, 0)\right)u_x^{l_2}(t) +\frac{1}{a}u_x^{l_2}(t)\ast\left(  \frac{g'_2(t)}{\alpha_2}- 2r'(l_2, t)\right) +O \left( \Delta t \right), $$
where
$$ u_x^{l_2}(t)={\frac {bt+1}{ \sqrt{\pi \,t}}}. 
$$
By using Maple software for example, the analytical approximations in short time limits  ($t\in [0, \Delta t] $), can be computed from above expressions as:
\begin{equation}
	\begin{array}{cccc}
		u^a(x,t)\simeq r(x,t)=\displaystyle\frac{1}{\sqrt{\pi t}}
		\int_{0}^{10}\varphi(\xi)\left[\exp\left( -\frac{(\xi-x)^2}{ t}\right) \right] d\xi\\
		\displaystyle=\frac{1}{4}\left( ( x-10) {\rm erf}\left( \frac{-10+x}{\sqrt{t}}\right) + (-2 x+10) {\rm erf}\left( \frac{-5+x}{\sqrt{t}}\right)+x\: {\rm erf}\left( \frac{x}{\sqrt{t}}\right)\right)\\
		\displaystyle+\frac{\sqrt{t}}{4\sqrt{\pi}}\left(-2\,{{\rm e}^{-{\frac { \left( x-5 \right) ^{2}}{t}}}}+{{\rm e}^{-{
					\frac { \left( x-10 \right) ^{2}}{t}}}}+{{\rm e}^{-{\frac {{x}^{2}}{t}
		}}} \right); 
		\label{eq52}
	\end{array}
\end{equation}
and
\begin{equation}
	u_x^a(0, t)=-u_x^a(10, t)=-\frac{1}{2}{\rm erf}\left( \frac{10}{\sqrt{t}}\right)+{\rm erf}\left( \frac{5}{\sqrt{t}}\right)+O \left( \Delta t \right)
	\label{eq53}
\end{equation}
where $ {\rm erf}(x) $ is the Error Function defined by:
$$\displaystyle{\rm erf}(x)=\frac{2}{\sqrt{\pi}}\int_0^x \mathrm{e}^{y^2}dy.$$
\begin{figure}[t]
	\centering
	\includegraphics[width=140mm]{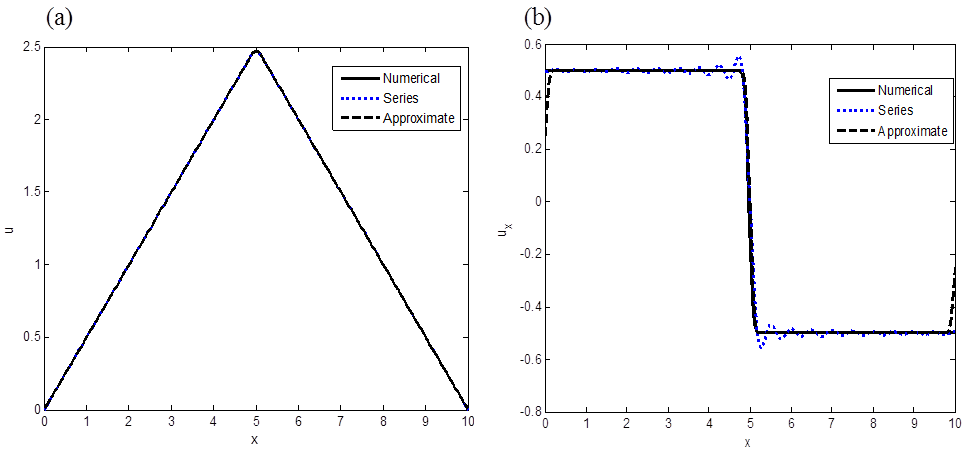}\hfill  
	\caption{(a) Solution $ u(x,t) $ at $ t=\Delta t=10^{-2}. $   (b) Partial derivative $ u_x(x,t)$ at $ t=\Delta t=10^{-2}$}
	\label{fig1}
\end{figure}
Compared to the function $ \varphi(x) $ characterizing the initial condition, it's remarkable that the derivatives coincide at $ l_1=0 $ as well as $ l_2=10 $ , that is :
\begin{equation}
	\displaystyle \lim_{t\rightarrow 0}u_x^a(0, t)=\varphi'(0)=\frac{1}{2},
	\label{eq54}	
\end{equation}
and
\begin{equation}
	\displaystyle \lim_{t\rightarrow 10}u_x^a(10, t)=\varphi'(10)=-\frac{1}{2}.	
	\label{eq55}	
\end{equation} 
In the limits of our knowledge, the truncation expansion (\ref{eq53}) was not available directly from the p-domain formula (\ref{eq51}). 

\begin{figure}[t]
	\centering
	\includegraphics[width=140mm]{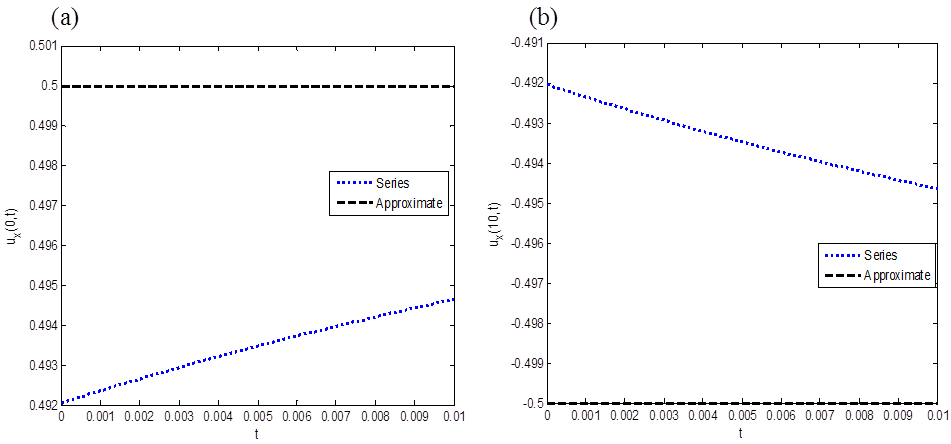}\hfill  
	\caption{(a) Evolution of the derivative $ u_x(0,t) $ for $t\in [0, \Delta t] $, $ \Delta t=10^{-2} $.   (b) Evolution of the derivative $ u_x(10,t) $ for $t\in [0, \Delta t] $, $ \Delta t=10^{-2} $}
	\label{fig2}
\end{figure}

Figure \ref{fig1} shows curves of the solution of the problem (\ref{eq44})-(\ref{eq46}) and of its derivative during the first time step $t\in [0, \Delta t] $ with $ \Delta t=10^{-2} $. The curves are obtained by using the three different methods. The numerical solution is computed by using toolbox functions of Matlab, while the series solution provided by the Fourier decomposition method is truncated at its first 20 terms. The approximate analytical solution is represented by the function $ u^a $ given by the above formula (\ref{eq52}). The three resulting curves of the solution $ u $ are almost identical as shown on figure \ref{fig1}(a). Their shapes, except around the peak point of abscissa $ x=5 $, seem to be very similar to that of the initial condition function (\ref{eq46}), since the time passed from $ t=0 $ to $ t=\Delta t=10^{-2} $ is still relatively small. On the figure \ref{fig1}(b), except at the domain boundaries, i.e. at $ x=0 $ and $ x=10 $, the derivative curves almost coincide for the approximate analytical and the numerical solutions. However, the curve for the series solution exhibits deviations in the form of small oscillations along the two former curves, especially around the peak point ($ x=5 $). This may suggest on the one hand, the potential consistency of the approximate analytical  solution as well as that of the numerical method to properly accounting for nonphysical phenomena as derivative jumps and the so-called Gibbs phenomenon arising in generalized series solutions. On the other hand, this may illustrate the pertinence of refined truncation expansions at the domain boundaries for the solution, that is for $ u^a(0,t) $ and $ u^a(10,t) $ which are keep null here, and for its derivative, i.e., $ u^a_x(0,t) $ and $ u^a_x(10,t) $ in the present case. Figures \ref{fig2}(a) and \ref{fig2}(b) respectively show evolution of the derivatives  $ u^a_x(0,t) $ and $ u^a_x(10,t) $ from formulae (\ref{eq53}). During the short time step $ [0, \Delta t] $, the time-varying derivative curves obtained from the approximate and series solutions can be compared. While the convergence of the series solution is poor especially at points close to the ends of the domain, the shapes shown by the curves of $ u^a_x(0,t) $ and $ u^a_x(10,t) $ correspond well to the similitude highlighted above with the initial function and confirmed by the calculations of the limits (\ref{eq54}) and (\ref{eq55}). Beyond its precision or consistency, another advantage of performing the approximate analytical solutions $ u^a_x(0,t) $ and $ u^a_x(10,t) $ during the short time limit $ [0, \Delta t] $, consists of the computational efficiency of these solutions. Indeed, at the first execution of the code source, the script running time of their procedure is far reduced compared to that of the series solution (about 20 times).  

Next, the exact solution in the Laplace domain can be extended to unbounded domains as infinite or semi-infinite intervals, for the space variable $ x $. Indeed when $ l_1=-\infty $ and $ l_2=+\infty $, one has $ U_x(l_1,p)=U_x(l_2,p)=0 $, and by taking the limits in equation (\ref{eq27}), the exact p-domain solution is reduced to $ U(x,p)=R(x,p) $, which corresponds to the following solution in the time domain:
$$
\begin{array}{ll} 
	u(x,t)=r(x,t)=\displaystyle\frac{\exp(-bt)}{2 a\sqrt{\pi t}}
	\int_{-\infty}^{+\infty}\varphi(\xi)\left[\exp\left( -\frac{(\xi-x)^2}{4 a^2 t}\right) \right] d\xi\\
	+\displaystyle\frac{1}{2 a\sqrt{\pi }}\int_{0}^{t}d\theta\int_{-\infty}^{+\infty}\frac{\exp(-b(t-\theta))}{\sqrt{(t-\theta)}}\exp\left( -\frac{(\xi-x)^2}{4 a^2 (t-\theta)}\right)f(\xi,\theta)d\xi.
\end{array} 
$$
In the case when $ b=0 $, the above solution is identical to the one reported by  \cite{Hen-partial} in their book at the section 6.8 titled: "The Heat Equation in an Infinite Region", where the Fourier decomposition method was used. About semi-infinite domains, let us chose $ l_1=0 $ and $ l_2=+\infty $ for example, then $ U_x(l_2,p)=0 $, and $ u(l_2,t) $ is to be considered as constant. The p-domain solution (\ref{eq27}) reduces to:
\begin{equation}
	U(x,p)=	\displaystyle \frac{1}{2}U(0,p)\exp\left(\frac{-x\sqrt{b+p}}{a} \right)-\frac{a U_x(0,p)}{2 \sqrt{b+p}}\exp\left(\frac{-x\sqrt{b+p}}{a} \right) + R(x,p), 
	\label{eq56}
\end{equation} 
where $R(x,p)={\mathcal L}\{r(x,t)\} $, and
\begin{equation}
	\begin{array}{ll} 
		r(x,t)=\displaystyle\frac{\exp(-bt)}{2 a\sqrt{\pi t}}
		\int_{0}^{+\infty}\varphi(\xi)\left[\exp\left( -\frac{(\xi-x)^2}{4 a^2 t}\right) \right] d\xi\\
		+\displaystyle\frac{1}{2 a\sqrt{\pi }}\int_{0}^{t}d\theta\int_{0}^{+\infty}\frac{\exp(-b(t-\theta))}{\sqrt{(t-\theta)}}\exp\left( -\frac{(\xi-x)^2}{4 a^2 (t-\theta)}\right)f(\xi,\theta)d\xi.
	\end{array} 
	\label{eq57}
\end{equation} 	
Setting as in the statement of the Problem a. of the chapter 8.1 in \cite{Lu-analytical}, $ b=0 $, $ \varphi(x)=t_0 =\mathrm{(constant)}$,  
$ u(0,t)=g_1(t)=t_a=\mathrm{(constant)}$ implying  $ U(0,p)=t_a/p $, it can be verified for $f(x,t)={w}/{c\gamma}=\mathrm{(constant)}$ that:
$$
Ux(0,p)=\displaystyle\frac{2\sqrt{p}}{a}R(0,p)-\frac{t_a}{a\sqrt{p}};
$$
and,
$$
R \left( x,p \right) =\displaystyle{\frac{t_0}{p}}+{\frac {w}{{p}^{2}c\gamma}}
-\,{\frac {t_0}{2p}{\exp\left( {-{\frac {x \sqrt{p}}{a}}}\right) }}-\,{
	\frac {w}{2{p}^{2}c\gamma}{\exp\left( {-{\frac {x \sqrt{p}}{a}}}\right) }}.
$$
This leads to reduce the exact p-domain solution (\ref{eq56}) into the form of:
$$
U \left( x,p \right) =\displaystyle{\frac{t_0}{p}}+{\frac {w}{{p}^{2}c\gamma}}
+\,{\frac {(t_a-t_0)}{p}{\exp\left( {-{\frac {x \sqrt{p}}{a}}}\right) }}-\,{
	\frac {w}{{p}^{2}c\gamma}{\exp\left( {-{\frac {x \sqrt{p}}{a}}}\right) }};
$$
which is identical to the solution (8.1.11) reported  in \cite{Lu-analytical}, provided that $ a $ is replaced by $ \sqrt{a} $ as specified in the statement of the problem. Note in the latter reference that the Laplace transform method was directly used to solve the problem, since the initial condition is specified as a constant function.

Finally, as a perspective for numerical simulation purposes, schemes to be used for deriving the complete analytical approximation of the solution $ u(x,t) $ on $ [l_1, l_2]\times[0, T] $  may be obtained. The procedure for calculating analytical approximations during the first time step $t\in [0, \Delta t]$ can be repeated for any short time step by updating the source term $ f(x,t) $, the functions $ \varphi(x) $, $ g_1(t) $ and $ g_2(t) $ that are related to the initial and boundary conditions. Also, the largest class to which the above functions must belong in order to obtain the exact solution in the Laplace domain may be studied more precisely. Moreover, computational models for one-dimensional Stefan problems as reviewed in \cite{Ja-Vu-Ve-Va-comparison} can be efficiently handled when using the truncation expansions $ u^a(l_1,t)$, $ u^a_x(l_1,t)$, $ u^a(l_2,t)$ and $ u^a_x(l_2,t)$ of the approximate analytical solution at the boundaries of the domain as expressed in section \ref{sec:5}. An example of using alike formulae for a specific problem of a spherically symmetric droplet evaporation can be seen in \cite{An-efficient}.

\section{Conclusion}
This study has permitted to calculate an explicit solution in the Laplace domain and analytical approximations in the earlier time-step to initial boundary value problems for the one-dimensional parabolic equation with constant coefficients. The problem is solved in its most general form, with the boundary conditions stated in a unified way on any  bounded generic interval of the real line. Compared to the classical or generalized series solutions that can be obtained for the same problem by using the Fourier decomposition method, the analytical approximations in short time limits are proven to be more consistent and sufficiently simple to improve computational efficiency in numerical schemes and simulations. Early time behaviors of heat or mass reaction-diffusion processes are of great interest in engineering and have a wide range of applications in fields like Computational Fluid Dynamics (CFD) and Nuclear Energy. In addition, the exact operational solution obtained for the problem can be extended to unbounded domains. This explicit solution can also be considered as an actual advance in the study of linear parabolic equations. Though they may be sought by means of Laplace inversion theorems, the number of exact inverse transforms from the Laplace domain into the time domain are limited in most tables or by using software facilities. However, inverse Laplace transforms can be accomplished numerically regardless of the complexity of the operational solutions. Therefore, accurate curves of the complete solution in the time domain can be obtained from the exact operational solution, when using a numerical inverse Laplace transform.

\bibliographystyle{unsrt}

\end{document}